%% file: ErrorMultiDelam.tex
\documentclass{article}

\usepackage{geometry}
\input{notations}
\usepackage{stmaryrd}
\usepackage{comment}

\title{On the control of the load increments for a proper description of multiple delamination in a domain decomposition framework.}

\author{O.~Allix, P.~Kerfriden, P.~Gosselet \vspace{5pt} \\ 
\large \textit{LMT-Cachan (ENS Cachan/CNRS/UPMC/PRES UniverSud Paris),} \\
\textit{61 av. du Pr\'esident Wilson, F-94230 Cachan, France}}

\date{January, 2010}

\begin{document}

\maketitle

In quasi-static nonlinear time-dependent analysis, the choice of the time discretization is a complex issue. The most basic strategy consists in determining a value of the load increment that ensures the convergence of the solution with respect to time on the base of preliminary simulations. In more advanced applications, the load increments can be controlled for instance by prescribing the number of iterations of the nonlinear resolution procedure, or by using an arc-length algorithm. These techniques usually introduce a parameter whose correct value is not easy to obtain. In this paper, an alternative procedure is proposed. It is based on the continuous control of the residual of the reference problem over time, whose measure is easy to interpret. This idea is applied in the framework of a multiscale domain decomposition strategy in order to perform 3D delamination analysis.

\section{INTRODUCTION}

The virtual testing of delamination is an objective widely spread among industrialists especially in the aeronautical field. To achieve it, two research thematics which have undergone large evolution during the last twenty years need to be put in conjunction: the pertinent modeling of composites and the efficient computation of structures.

Indeed, there have been many advances toward a better understanding of the mechanics of laminated composites and of damage mechanisms. Two kinds of modeling have proved their validity: microscale and mesoscale models. Microscale models are strongly connected to the physics of the material and thus provide a reliable framework for simulation \cite{ladeveze05,lubineau09}. 
Unfortunately, the computation of models defined at the micro scale require such a fine discretization that only small test specimens can be simulated, structural computations being out of reach even on recent hardware. Meso-models \cite{schellekens94,allix98,ladeveze02c,deborst06} are defined at a scale which enables both the introduction of physics-based ingredients and the simulation of small industrial structures. They indeed most often rely on the definition of two meso-constituents, the ply (3D entity) and the interface (2D entity), which are modeled using continuum (damage) mechanics, their behavior being obtainable from the homogenization of micro-models \cite{lubineau07}. Anyhow, for reliable simulation, discretizations still need to be fine (in order, for instance, to represent correctly the gradients of stresses due to edge effects which are responsible for the initiation of many degradations) and associated systems thus remain very large (in terms of number of degrees of freedom) and strongly nonlinear (with potential instabilities). 

As a first approach of the reliable simulation of quasi-static simulations of the delamination in composite structures, we chose in \cite{kerfriden09} to neglect the effect of deterioration within the plies and to lump the degradations in the interfaces. We thus retained the mesomodel presented in \cite{allix92} where the delamination ability is localized in the interfaces and handled through a cohesive behavior. The space discretization is considered sufficiently fine to represent accurately any evolution of multiple delamination cracks (sufficient number of Gauss points in the length of the process zone \cite{schellekens93,alfano01,deborst06}).
At each time step of an incremental time discretization scheme, the associated large nonlinear system is solved using a three-scale domain decomposition strategy. Based on the mixed LaTIn-based domain decomposition method \cite{ladeveze00}, this strategy has been given high numerical efficiency by adapting various ideas from the work of \cite{mandel93,moes99,pebrel08} to the computation of delamination. Three-dimensional simulations of the delamination in realistic structural components have been performed on parallel computers without the need to perform local space refinement. 

Though, a complex issue arises when choosing the load increments: the solution to softening quasi-static problems depends on the time discretization scheme parameters (non-uniqueness of the solution and possible bifurcation paths). This remark brings us to the field of the validation. In the literature, numerous error indicators have been developed to control \textit{a posteriori} the global error introduced in finite element schemes for linear problems \cite{babuska78,zienkiewicz87,ladeveze91}. These indicators have been extended to the validation of nonlinear time-dependent problems \cite{ladeveze86b,bass87,Johnson92}. One of the most advanced criterion is the so-called error in the constitutive law \cite{ladeveze86b}. A solution to the nonlinear evolution problem being computed using a FE scheme and a classical time integration procedure, one constructs a solution which satisfies the kinematic and static admissibilities, and lump the residual of the nonlinear evolution problem equations in the constitutive laws. A measure of this residual permits to control at the same time the discretization error in space, in time and the error introduced by the iterative solution procedures \cite{ladeveze99c}. This idea has been formalized in \cite{ladeveze99} for materials described using internal variables. The state equations are satisfied by the reconstructed solutions, the measure of the non-verification of the evolution laws permits to derive a strict upper bound to the solution error. Though, this new admissible solution is not easily constructed in the case of softening behaviors. Specific developments in \cite{ladeveze99b} meant to tackle this difficulty, and the resulting procedure is used in \cite{ladeveze99c} to derive an adaptive refinement procedure in space and time. Note that, at the present time, a link between the error in the constitutive law and the error in the solution is still to be established in the case of softening materials.

The aim of the work presented in this paper is double. The first is to define a comprehensive time discretization error indicator inspired from the work of \cite{ladeveze86b,ladeveze99c} for delamination analysis and to ensure that its computation and use is numerically efficient within the LaTIn-based domain decomposition strategy. Our second goal is to use the developed indicator to control on the fly the load increment in quasi-static analysis in order to ensure the convergence of the computed solution.

The paper is organized as follows. The reference delamination problem is presented in Section \ref{sec:reference_problem}. The dependency of the solution to this problem on the time discretization scheme is  demonstrated. In the following section, we present a time-dependent error indicator based on the error in the constitutive law and computed with respect to a continuous solution in time, constructed by interpolation over each time step. Although very general, this indicator is not directly suitable for the LaTIn-based multiscale strategy used to perform the nonlinear resolutions. The main features of this strategy are recalled in Section \ref{sec:latin}. We focus in particular on the indicator based on the error in the constitutive law used to estimate the convergence of the iterative procedure. In Section \ref{sec:time_control}, this convergence indicator is associated to the previous developments to derive an alternative and cheap time discretization error indicator, which is the basis for the development of an automatic time-step-control procedure. At last, this technique is validated on multiscale and parallel delamination simulations in Section \ref{sec:results}. Two different problems are assessed: a simple and stable problem in which the time increments correspond to the increases in the prescribed load, and a more complex and unstable problem, solved using an arc-length procedure, in which the time increments correspond to the value of the arc-lengths.

\section{THE REFERENCE PROBLEM AND ITS DISCRETIZATION IN TIME}
\label{sec:reference_problem}

\input{chap_0}

\section{A TIME DISCRETIZATION ERROR INDICATOR}
\label{sec:criterion}

\input{chap_1}

\section{THE NONLINEAR RESOLUTION STRATEGY}
\label{sec:latin}

\input{chap_2}

\section{AN AUTOMATIC PROCEDURE TO CONTROL THE LOAD INCREMENTS}
\label{sec:time_control}

\input{chap_3}

\section{VALIDATION OF THE STRATEGY}
\label{sec:results}

\input{chap_4}

\section{CONCLUSION}

In this paper, we presented a strategy to adapt automatically the time increment in quasi-static delamination problems to the very sharp non-linearities which are involved. This strategy is based a continuous monitoring of the residual of the reference problem equations with respect to time. This has been achieved by calculating the error in the constitutive law on admissible solutions interpolated over each time steps, which enables to define a time discretization error criterion evaluating the relevancy of the nonlinear computations performed at each time increment. Based on this indicator, a simple procedure to control the time step has been derived. The main parameter of this technique is easy to obtain as it only requires to perform time-independent benchmark tests prior to the delamination simulations. The validity of this procedure has been demonstrated on delamination problems undergoing global instabilities. 

Our current interest being to perform buckling-driven delamination analysis, the validity of this strategy shall be verified, in the future, on computations involving geometrical nonlinearities.




\bibliographystyle{plain}
\bibliography{template}

\end{document}

%% file: chap_0.tex
\subsection{Reference problem at a given time of the analysis}



The delamination simulation is performed under the assumptions of quasi-static, isothermal evolution over time and small perturbations.

The laminate structure $\mathbf{E}$ occupying Domain $\Omega$ is made out of $N_P$ adjacent plies occupying Domains $(\Omega_P)_{P \in \llbracket 1, \ N_P\rrbracket}$ (of boundaries $(\partial \Omega_P)_{P \in \llbracket 1, \ N_P\rrbracket}$), separated by $(N_P-1)$ cohesive interfaces $(I_{P})_{P \in \llbracket 1, \ N_P-1\rrbracket}$ and  (see Figure (\ref{fig:decomp_sst_interfaces}), Page \pageref{fig:decomp_sst_interfaces}).  An external traction field $\V{F}_d$ (respectively a displacement field $\V{U}_d$) is applied to the structure on Part $\partial \Omega_f$ (respectively $ \partial \Omega_u$) of the boundary $\partial \Omega$ of Domain $\Omega$. The volume force is denoted $\V{f}_d$. Let $\V{u}_P$ be the displacement field, $\M{\sigma}_P$  the Cauchy stress tensor and $\M{\epsilon}_P$ the symmetric part of the displacement gradient in Ply $P$. 

At every time $t \in [0 \ T]$ of the analysis, the reference non-linear equilibrium problem reads: 

\vspace{5pt}
\textit{Find $\displaystyle s_{ref} = (s_P)_{P \in \llbracket 1, \ N_P\rrbracket)} $, where $\displaystyle s_P= (\V{u}_P, \M{\sigma}_P)$, which satisfies the following equations:}
\begin{itemize}
\item Kinematic admissibility on $\partial \Omega_u$:
\begin{equation}
\label{eq:kinematic_admissibility}
{\V{u}_P}_{| \partial \Omega_u} = \V{U}_d
\end{equation}
\item Global equilibrium of Structure $\mathbf{E}$: $\forall \displaystyle {({\V{u}_P}^{\star})}_{P \in \llbracket 1, \ N_P\rrbracket}$
\begin{equation}
\label{eq:static_admissibility}
\begin{array}{ll}
\displaystyle \sum_P   \int_{\Omega_P}   \operatorname{Tr} \left( \M{\sigma}_P \M{\epsilon}({\V{u}_P}^{\star}) \right) d \Omega   
 & \displaystyle - \sum_P \int_{\Omega_P} \V{f}_d . {\V{u}_P}^{\star} \ d\Omega - \sum_P \int_{\partial \Omega_P \cap \partial \Omega_f} \V{F}_d . {\V{u}_P}^{\star} \ d\Gamma
\\
& \displaystyle + \sum_{P} \int_{I_{P}} \M{\sigma}_P \V{n}_P . {\V{[ u ]}_P}^{\star} \ d\Gamma = 0
\end{array}
\end{equation}
where $\V{[u]}_P$ is the jump of displacement of Interface $I_{P}$: $\V{[u]}_P = \V{u}_{P+1}-\V{u}_{P}$ and   $\V{n}_P$ is the outer normal to the boundary $\partial \Omega_P$.
\item Linear orthotropic behavior of the plies:
\begin{equation}
\label{eq:constitutive_ply}
\quad \M{\sigma}_P = \mathcal{K} \, \M{\epsilon}(\underline{u}_P)
\end{equation}
\label{eq:constitutive_interface}
\item Constitutive law of the cohesive interfaces, local on any interface $I_{P}$. The elastic damageable law proposed in \cite{allix98} is described using continuum damage mechanics. Three internal variables $(d_i)_{i \in \llbracket 1, \ 3\rrbracket}$ (one for each delamination mode: traction along $\V{n}_P$ and shear along $\V{t}_1$ and $\V{t}_2$ on Figure (\ref{fig:interface_composite})), ranging from 0 to 1 are introduced in the surface strain energy $e_d$ in order to take into account the irreversible damage mechanisms.
\begin{figure}[ht]
       \centering
       \includegraphics[width=0.65 \linewidth]{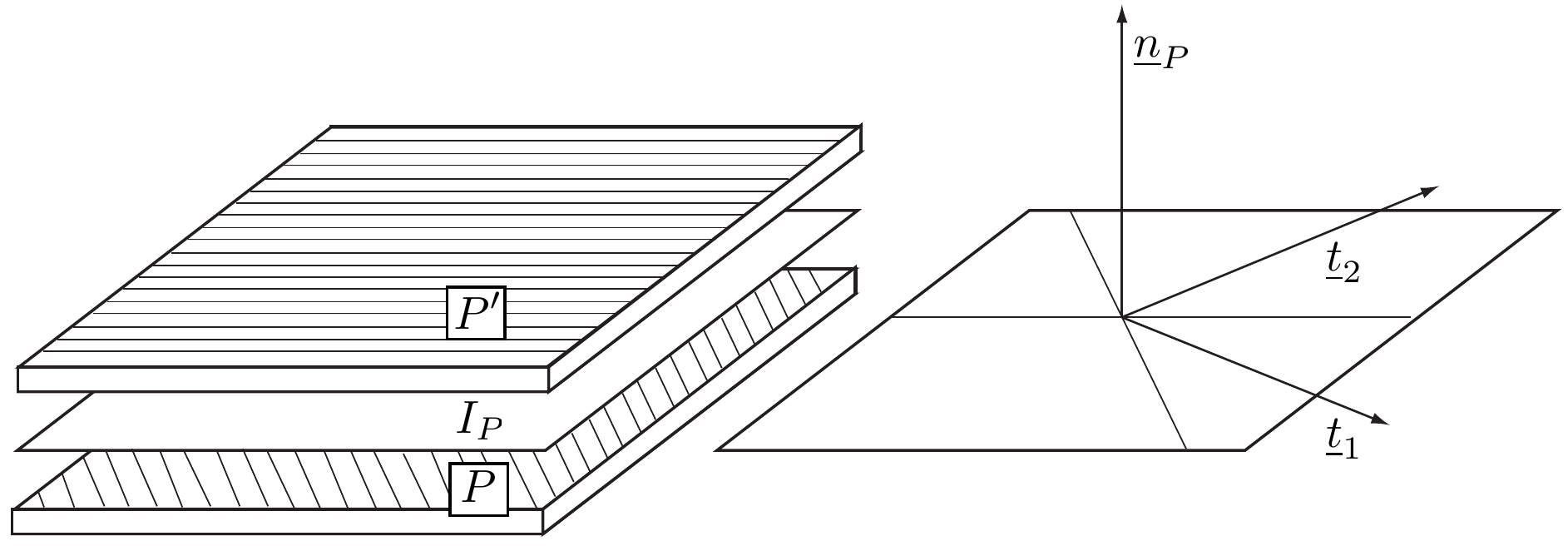}
       \caption{The mesomodel entities}
       \label{fig:interface_composite}
\end{figure}
\begin{itemize}
\item Two state equations are derived from the expression of the free energy. The first one establishes a relation between the dual interface unknown $\M{\sigma}_P \V{n}_P$, and the primal interface unknown $\V{[u]}_P$:
\begin{equation}
\label{eq:elast_endom}
\displaystyle \M{\sigma}_P.\V{n}_P = \frac{\partial e_d}{\partial \V{[u]}_P} \quad \textrm{which gives} \quad \M{\sigma}_P.\V{n}_P = \M{K}_{P} \left( \left(\V{[u]}_P \right)_{|\tau \in [0 \  t]}  \right) \V{[u]}_P 
\end{equation}
where, in the basis $(\V{n}_P,\V{t}_1,\V{t}_2)$, $h_+$ being the positive indicator function:
\begin{equation}
\nonumber
\M{K}_{P}  \left( \left(\V{[u]}_P \right)_{|\tau \in [0 \  t]}  \right) = \left( \begin{array}{ccc}
\displaystyle \left(1-h_+ ( \V{[u]}_P.\V{n}_P ) d_3 \right) k_n^0  & 0 & 0 \\
0 & \displaystyle (1-d_1) k_t^0 & 0 \\
0 & 0 & \displaystyle (1-d_2) k_t^0 
\end{array}
\right)
\end{equation}
The second state equation links the thermodynamic forces $(Y_i)_{i \in \llbracket 1, \ 3 \rrbracket}$ to the primal interface unknown:
\begin{equation}
Y_i = - \frac{\partial e_d}{\partial d_i}
\qquad \textrm{where} \qquad
\left\{ \begin{array}{ccl}
Y_1 & = & \displaystyle \frac{1}{2} \, k_t^0 \, \left(\V{[u]}_P.\V{t}_1 \right)^2 \\
Y_2 & = & \displaystyle \frac{1}{2} \, k_t^0 \, \left(\V{[u]}_P.\V{t}_2] \right)^2 \\
Y_3 & = & \displaystyle \frac{1}{2} \, k_n^0 \, \left( h_+(\V{[u]}_{P}.\V{n}_P ) \right)^2
\end{array} \right. \end{equation}
\item The evolution laws are:
\begin{equation}
\label{eq:evolu_law}
\begin{array}{l}
\displaystyle d_1 = d_2 = d_3 = \min \{ 1 ,w(Y) \} \\ \textrm{where}
\ \left\{ \begin{array}{ll} 
w(Y) = \frac{n}{n+1} \left( \frac{Y}{Y_c} \right)^n \\
Y =  \operatorname{max}_{(\tau \leq t)} \left( {Y_3}_{|\tau}^\alpha + \gamma_1 {Y_1}_{|\tau}^\alpha + \gamma_2 {Y_2}_{|\tau}^\alpha \right)^{\frac{1}{\alpha}}
\end{array} \right.
\end{array}
\end{equation}
\end{itemize}
Further details on this cohesive zone model and identification issues can be found in \cite{allix98}.

The dissipated energy $E_{dissi}$ will be used in this paper as a global measure of the delaminated area of the cohesive interfaces:
\begin{equation}
E_{dissi} =    \sum_{P} \int_{I_{P}} \int_{0}^t \left( \sum_{i=1}^3 \, Y_i \, \dot{d} \right) \ dt \, d \Gamma  = \sum_{P}  \int_{I_{P}} A \, d \ d \Gamma 
\end{equation} 
where A is a scalar which only depends on the parameters of the interface model.
\end{itemize}

In the following developments, the investigations are restricted to simulations under prescribed forces and displacements following a unique load function of time. In this context, the volume force will be assumed negligible. These assumptions are not mandatory to make use of the work presented in this paper, but simplify the construction of a continuous solution over time (Section \ref{sec:criterion}).

\subsection{Time discretization scheme}

An incremental procedure is used to solve the problem over time. It consists in discretizing the time of the analysis $[0 \ T]$ in $N$ intervals $[t_n \ t_{n+1}]_{n \in  \llbracket 0, \ N-1\rrbracket}$. Successive nonlinear problems are solved at each computation time $(t_n)_{n \in \llbracket 0, \ N \rrbracket}$. 

Hence, a solution to the discretized problem in time is a set of $N+1$ solutions satisfying the reference problem equations, the time dependency in the constitutive laws being discretized. More precisely, at Computation time $t_{n+1}$, the discretization of Equations \eqref{eq:elast_endom} and \eqref{eq:evolu_law} reads:
\begin{equation}
\label{eq:elast_endom_discr}
\displaystyle \M{\sigma}_P.\V{n}_P = \M{K}_{P} \left( \left(\V{[u]}_P \right)_{| t \in [ t_0, \  t_{n+1} ]} \right) \V{[u]}_P 
\end{equation}

In general, the time discretization is chosen so that within each interval $[t_n \ t_{n+1}]_{n \in  \llbracket 0, \ N-1\rrbracket}$, the evolution of the prescribed load is monotonic, which will also be assumed in the following. 



\subsection{Influence of the time increments on the solution to the discretized delamination problem}

The solution to the discretized reference problem reached at time $T$ strongly depends on the time increments for two reasons: 
\begin{itemize}
\item the discretized cohesive law (Equation \eqref{eq:elast_endom_discr}) depends on the discrete history of the interface variables. Hence, the residual stiffness of the cohesive interfaces depends on the time increments. This phenomenon is illustrated in the next section.
\item structural problems involving softening materials may be unstable and may have multiple solutions. In those cases,  the solution paths depend on both the algorithm used at each computation time step and the initialization of this algorithm (\textit{i.e.}: the previous converged solution). The resulting dependency on the time increments will be demonstrated in the last section of this paper.
\end{itemize}

\paragraph{DCB (double cantilever beam) test case}

\begin{figure}[htb]
       \centering
       \includegraphics[width=0.9 \linewidth]{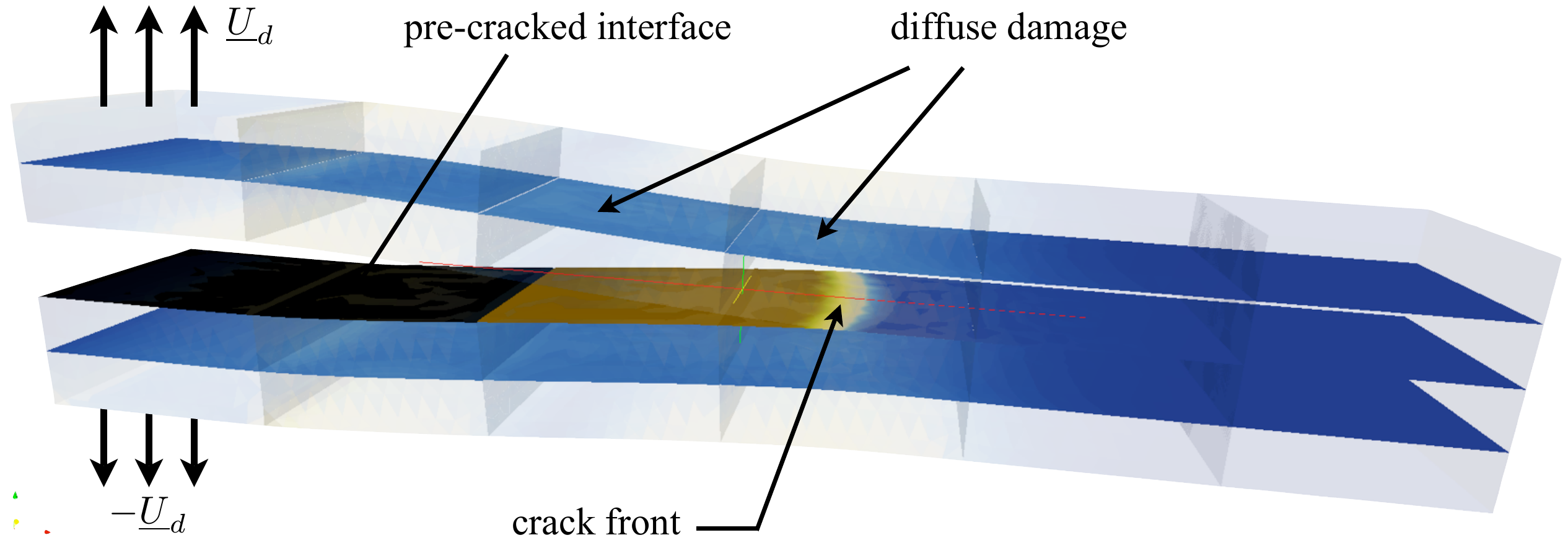}
       \caption{Definition of the four-ply DCB problem}
       \label{fig:DCB}
\end{figure}

\begin{figure}[ht]
       \centering
       \includegraphics[width=0.6 \linewidth]{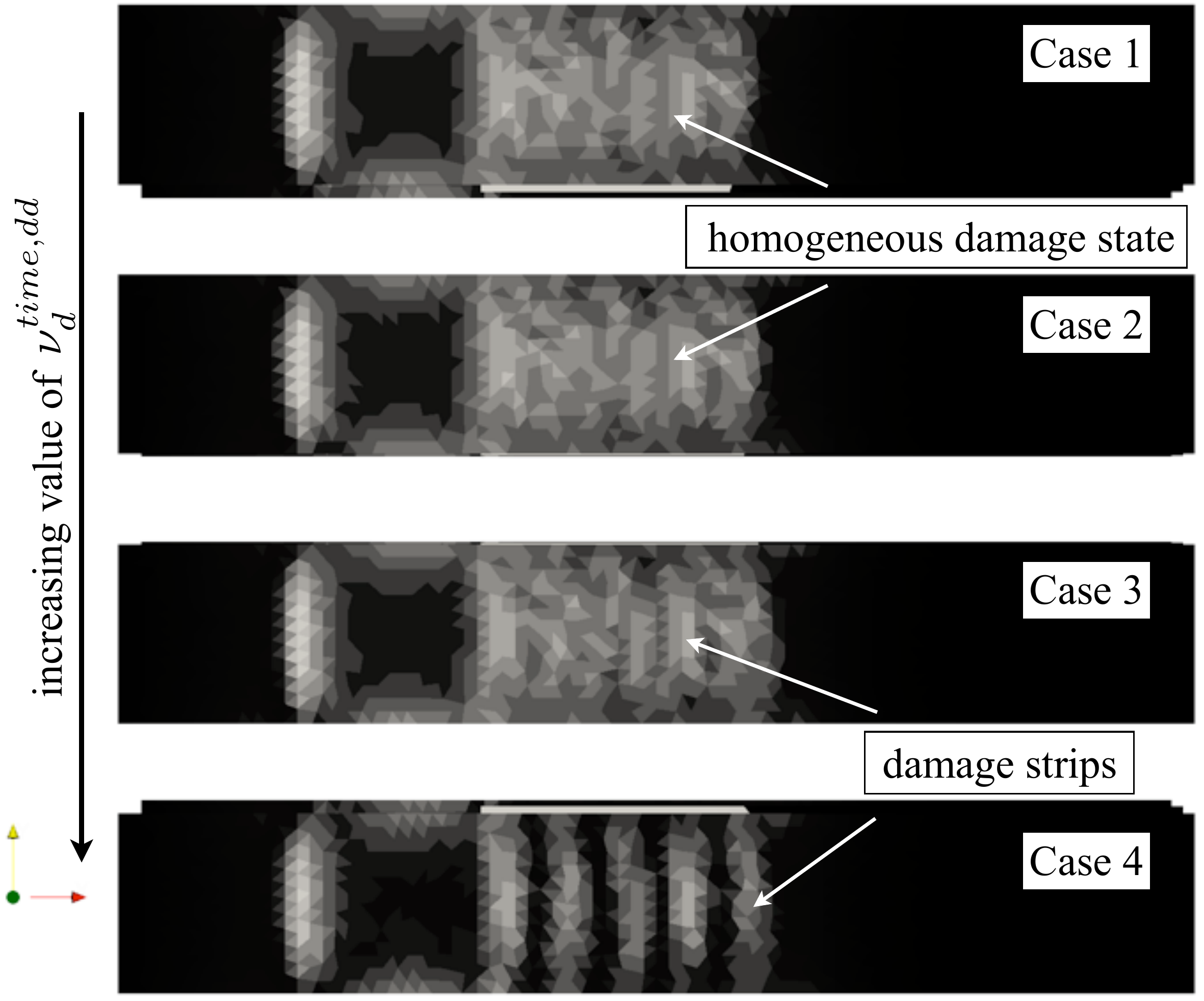}
       \caption{Influence of the prescribed value $\nu_d^{time,dd}$ on the damage state in the upper cohesive interface of the  DCB problem}
       \label{fig:DCB_time}
\end{figure}

The laminate structure that we consider is made out of four isotropic plies (Figure \ref{fig:DCB_time}). One part of the median cohesive interface is replaced by a contact interface in order to simulate an initial crack in the structure. Displacements are prescribed for the crack to propagate in a stable manner. The final prescribed displacement is set to a predefined value, which fixes the propagation length. The initial stiffness of the cohesive interfaces is obtained by integrating the Young and shear moduli of the matrix in the ``thickness'' of the interface ($1/10$ the thickness of the plies) \cite{allix98}.

The solution is not unique and depends on the load increments. Figure \ref{fig:DCB_time} presents the damage state in the upper cohesive interface, four different time discretizations being applied (these results will be fully detailed later on, for the values of the successive load increments are obtained by the adaptive time step procedure described in Section \ref{sec:time_control}). $\nu_{d}^{time,dd}$ is the criterion driving the time discretization (the largest $\nu_{d}^{time,dd}$, the coarser the discretization). In cases 1 and 2, the number of time increment used in the propagation phase of the analysis are, respectively, 69 and 21. The differences in the damage state of the interfaces are not significant, the evolution of the crack front being sufficiently slow to capture the effects of the stress concentrations. Hence, both these solutions are converged with respect to the time. In case 3, obtained with 9 coarse time increments, the solution is slightly different from the previous reference cases. Finally, in case 4, using only 5 time increments to describe the propagation of the crack clearly leads to the appearance of damage strips in the upper and lower interfaces. This is due to the effect of the stress concentration at the tip of the crack which propagates in a discrete manner with respect to time.

%% file: chap_1.tex
We suppose that two consecutive solutions to the reference problem, $s_n$ at Time $t_n$ and $s_{n+1}$ at Time $t_{n+1}$, have been computed using a nonlinear resolution strategy. The aim is to evaluate the relevancy of the solution computed at Time $t_{n+1}$, the continuous evolution of the structure over the current time step $[ t_n \ t_{n+1} ]$ being \textit{a priori} unknown. We propose to construct an interpolated solution over the time step in order to monitor the residual of the nonlinear reference problem continuously.

\subsection{Interpolation of the kinematic and static fields over a time step}

Let us prescribe the continuous evolution of the prescribed boundary values over the time step: 
\begin{equation}
\forall \, {\bar t} \in [t_n \ t_{n+1}], \qquad
\left \{\begin{array}{l}
\displaystyle \forall \, M \in \partial \Omega_f, \qquad {\V{F}_d}_{|  \bar t} =  \alpha({\bar t}) \, {\V{F}_d}_{| t_n}  + (1-\alpha({\bar t})) \, {\V{F}_d}_{| t_{n+1}} \\
\displaystyle  \forall \, M \in \partial \Omega_u, \qquad {\V{U}_d}_{|  \bar t} =  \alpha({\bar t}) \, {\V{U}_d}_{| t_n}  + (1-\alpha({\bar t})) \, {\V{U}_d}_{| t_{n+1}}
\end{array} \right.
\end{equation}
where the function $\alpha(\bar {t})$ is the restriction of the load function over $[t_n \ t_{n+1}]$. In the case of a linear evolution (which will be the case in our applications), it simply reads:
\begin{equation}
\forall \, {\bar t} \in [t_n \ t_{n+1}], \qquad \alpha(\bar {t}) = \frac{{\bar t}-t_n}{t_{n+1}-t_n}
\end{equation}

\begin{figure}[ht]
       \centering
       \includegraphics[width=0.53 \linewidth]{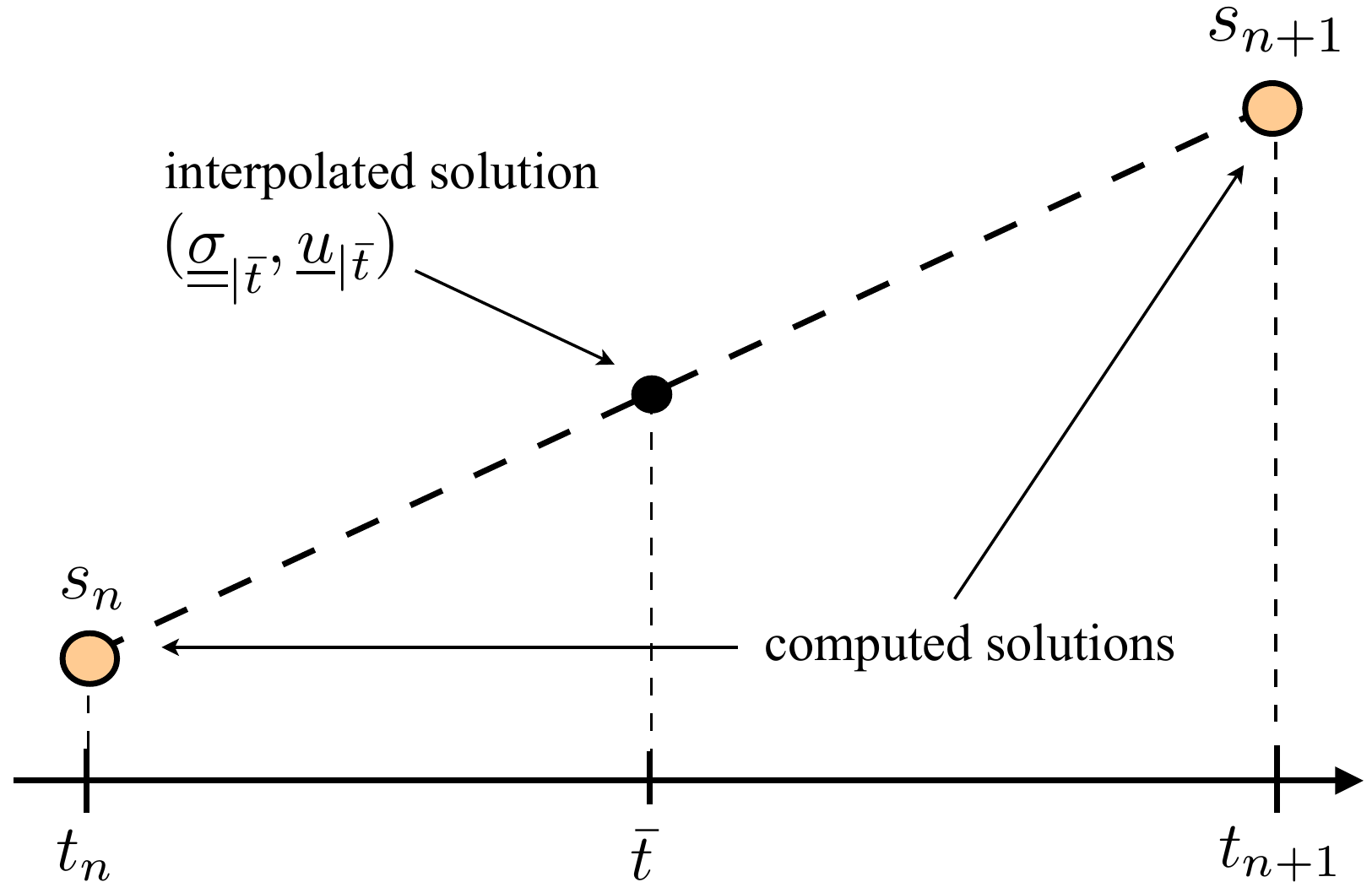}
       \caption{Schematic representation of the interpolation performed over each time step}
       \label{fig:time_criterion}
\end{figure}

The evolution of the kinematic and static fields over the current time is assumed to follow the evolution of the prescribed loading (see Figure (\ref{fig:time_criterion})), which writes:
\begin{equation}
\forall \, {\bar t} \in [t_n \ t_{n+1}], \, \forall P \in \llbracket 1, \ N_P \rrbracket, \qquad
\left\{\begin{array}{l}
\displaystyle {\V{u}_P}_{|  \bar t} = \alpha({\bar t}) \, {\V{u}_P}_{| t_n}  + (1-\alpha({\bar t})) \, {\V{u}_P}_{| t_{n+1}} \\
\displaystyle {\M{\sigma}_P}_{|  \bar t} = \alpha({\bar t}) \, {\M{\sigma}_P}_{| t_n}  + (1-\alpha({\bar t})) \, {\M{\sigma}_P}_{| t_{n+1}} 
\end{array} \right.
\end{equation}

$s_n$ and $s_{n+1}$ are two solutions of the reference problem. In particular, they satisfy the following set of linear equations:
\begin{itemize} 
\item kinematic admissibility, Equation \eqref{eq:kinematic_admissibility}
\item static admissibility, Equation \eqref{eq:static_admissibility}, the volume force being assumed negligible.
\item constitutive law of the plies, Equation \eqref{eq:constitutive_ply}
\end{itemize}
As a consequence, the interpolated kinematic and static fields over the current time step also satisfy this set of linear equations. Hence, the residual of the reference problem at any time $\bar t \in [t_n \ t_{n+1}]$ is the residual of the constitutive law of the cohesive interfaces, which remains the only non-satisfied equation.

\subsection{Evolution of the damage variables over the current time step}

At any intermediate time $\bar{t} \in [t_n \ t_{n+1}]$, the internal variables are calculated with respect to the continuous history of the interpolated displacement field on Time interval $[0 \ \bar{t}]$. Let us define a new stress field $\M{\widehat{\sigma}}$ which satisfies the nonlinear constitutive law of the interfaces:
\begin{equation}
\begin{array}{r}
\displaystyle \forall \, {\bar t} \in [t_n \ t_{n+1}], \, \forall \, P \in \llbracket 1, \ N_P-1\rrbracket, \ \textrm{on} \ I_{P}, \\
\displaystyle {\M{\widehat{\sigma}}_P}_{| \bar t}.\V{n}_P = \M{K}_{P} \left( \left(\V{[u]}_P \right)_{|\tau \in [0 \ \bar t]}  \right) \ \V{[u]}_{|\bar t}
\end{array}
\end{equation}

Alternatively, one can update the damage variables with respect to the interpolated stress field, and define a jump of displacement field $\V{\widehat{[u]}}$ 
satisfying the constitutive law of the cohesive interfaces.

The damage variables initially computed at time $t_{n+1}$ by the nonlinear resolution strategy are discarded. Indeed, they may differ from the ones obtained at time $t_{n+1}$ by  the continuous construction over $[t_n \ t_{n+1}]$, for solution $s_{n+1}$ only satisfies the discretized cohesive law \eqref{eq:elast_endom_discr}. The residual of the reference problem equations at Time $t_{n+1}$ obtained when updating the damage variables can be reduced by lowering the time increment $\Delta t = t_{n+1}-t_n$ and performing new nonlinear resolutions at Time $t_{n+1}$, which will be detailed in Section \ref{sec:time_control}.

\subsection{Definition of the time discretization error indicator}
 
A measure  $\nu^{interp}$ (``\it{interp}'' \normalfont stands for ``interpolation'') of the residual of the reference problem equations at any time $\bar t \in [t_n \ t_{n+1}]$ can be obtained by summing the local contributions to the error in the nonlinear constitutive laws:
\begin{equation}
\label{eq:nu}
\nu^{interp}_{| \bar t} = \sum_{ P } \frac{\parallel \left({\M{\sigma}_P}_{| \bar t} - {\M{\widehat{\sigma}}_P}_{| \bar t} \right) \V{n}_P \parallel_{I_{P} } }{\parallel \left({\M{\sigma}_P}_{| \bar t} + {\M{\widehat{\sigma}}_P}_{| \bar t} \right) \V{n}_P \parallel_{I_{P} } }
\qquad \textrm{where} \quad \parallel \ x \ \parallel_{I_{P} } = \int_{I_{P}} x^T \ x \ d \Gamma
\end{equation}
Or alternatively if the history is updated with respect to the interpolated stress field, 
\begin{equation}
\label{eq:nutilde}
\widetilde{\nu}^{interp}_{| \bar t} = \sum_{ P } \frac{\parallel \left({\V{[u]}_P}_{| \bar t} - {\V{\widehat{[u]}}_P}_{| \bar t} \right)  \parallel_{I_{P}}}{\parallel \left({\V{[u]}_P}_{| \bar t} + {\V{\widehat{[u]}}_P}_{| \bar t} \right) \parallel_{I_{P}}}
\end{equation}
\begin{figure}[ht]
       \centering
       \includegraphics[width=0.59 \linewidth]{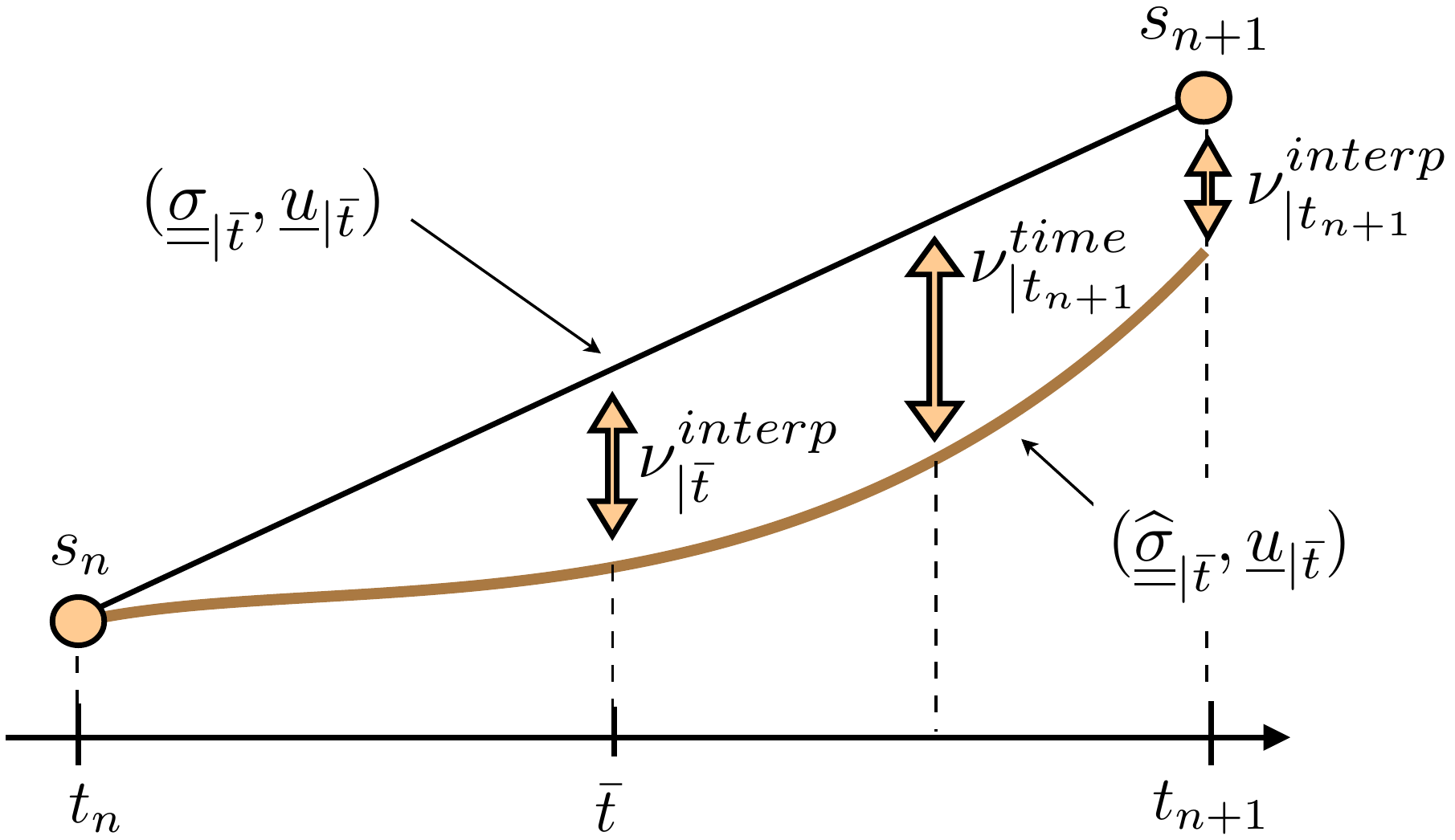}
       \caption{Schematic representation of the time discretization error indicator}
       \label{fig:time_criterion2}
\end{figure}

The time discretization error indicator at Time $t_{n+1}$ is defined as the maximum value of the previous measure over  $[t_n \ t_{n+1}]$ (see Figure (\ref{fig:time_criterion2})), which reads:
\begin{equation}
\nu^{time}_{| t_{n+1}}  = \max_{\bar t \in [t_n \ t_{n+1}]} \nu^{interp}_{| \bar t} \qquad \textrm{or alternatively} \qquad \widetilde{\nu}^{time}_{| t_{n+1}} = \max_{\bar t \in [t_n \ t_{n+1}]} \widetilde{\nu}^{interp}_{| \bar t}
\end{equation}

The concept introduced here finds its roots in the work of \cite{ladeveze86b,gallimard96}, in which the sum over time of the product of Criteria \eqref{eq:nu} and \eqref{eq:nutilde} is used to measure the error in the constitutive law due to both space and time discretization for materials satisfying Drucker's stability equality. Three main differences should be outlined here:
\begin{itemize}
\item In the case of softening materials, Drucker's stability equality is not satisfied. The mathematical properties which result from the definition of the Drucker's law-based criterion do not apply. Hence, making use of this criterion is not relevant. In addition, computing $\widetilde{\nu}^{time}$ requires the monotony of the interface behavior (uniqueness of the admissible displacement jump for any arbitrary local stress state). In the following developments, we will use Criterion $\nu^{time}$ to measure the residual of the reference problem equations over the current time step.
\item Our final goal being to provide an algorithm to control "on-the-fly" the time increments, $\nu^{time}$ is not a norm over the whole time of the analysis, but it instead is evaluated locally over each time increment.
\item To be consistent with \cite{ladeveze86b,gallimard96} the field ${\M{\widehat{\sigma}}_P}_{| \bar t}$ should also be reconstructed with respect to the space variables so that it satisfies exactly the static admissibility condition \eqref{eq:static_admissibility}. In this work we focus on the time discretization and so we content ourselves with a weak (discrete) static admissibility. At Times $t_n$ and $t_{n+1}$ solution fields satisfy the constitutive law of the plies \eqref{eq:constitutive_ply}, the kinematic admissibility and the static admissibility ``in the finite element sense''. Thus Criterion $\nu^{time}$ (which is introduced without reference to space discretization) only accounts for the error due to time discretization.
\end{itemize}

\subsection{Practical considerations}

\subsubsection{Sub-intervals}
In practice, $\nu^{interp}$ is computed at a given set of intermediate times within the current time step. $ [t_n \ t_{n+1}] $ is subdivided into $N_s$ subintervals $[\bar t_i \ \bar t_{i+1}]_{i \in \llbracket 0, \ N_s-1\rrbracket}$, the time discretization error indicator ${\nu_{|t_{n+1}}^{time}}$ being computed as:
\begin{equation}
{\nu_{|t_{n+1}}^{time}} = \max_{i \in  \llbracket  0, \ N_s  \rrbracket} \nu^{interp}_{|  \bar  t_i}
\end{equation}

\subsubsection{Error in the cohesive law}
Computing $\nu^{time}$ requires to extract the transverse constraints $(\M{\sigma}_P.\V{n}_P)_{P \in \llbracket 1, \ N_P\rrbracket}$ which is not directly available in finite element codes. Usually, cohesive interface elements are used to overcome this problem. Classical incremental Newton solvers can then be used to solve the delamination problem at each computation time $(t_n)_{n \in \llbracket 0, \ N \rrbracket}$. The technique to control the time increment that we propose in Section \ref{sec:toto} can directly be applied to such approaches.

We focus on the insertion of the control technique within the framework described in \cite{kerfriden09}. The principle is to use an incremental LaTIn-based domain decomposition strategy \cite{ladeveze03b} to efficiently solve (in parallel) the delamination problem at each computation time. In this case, the cohesive behavior is directly described as a nonlinear joint between substructures. The mixed description of the interface behavior makes the transverse constraints available naturally. As it shall be detailed in Section \ref{sec:time_control}, the time discretization error indicator can be defined as a time-dependent version of the convergence indicator used to stop the iterations of the LaTIn algorithm.


%% file: chap_2.tex
We propose an overview of the domain decomposition strategy used to perform the successive nonlinear resolutions of the delamination analysis, first in the stable case, then in the unstable case, where it is combined with an arc-length procedure. We focus in a second time on the development of a convergence indicator based on the error in the constitutive law \cite{ladeveze86b} to stop both of these iterative solvers. Further details concerning the multiscale and parallel computing aspects can be found in \cite{kerfriden09}.

\subsection{Substructured formulation of the reference problem}

\begin{figure}[ht]
       \centering
       \includegraphics[width=0.85 \linewidth]{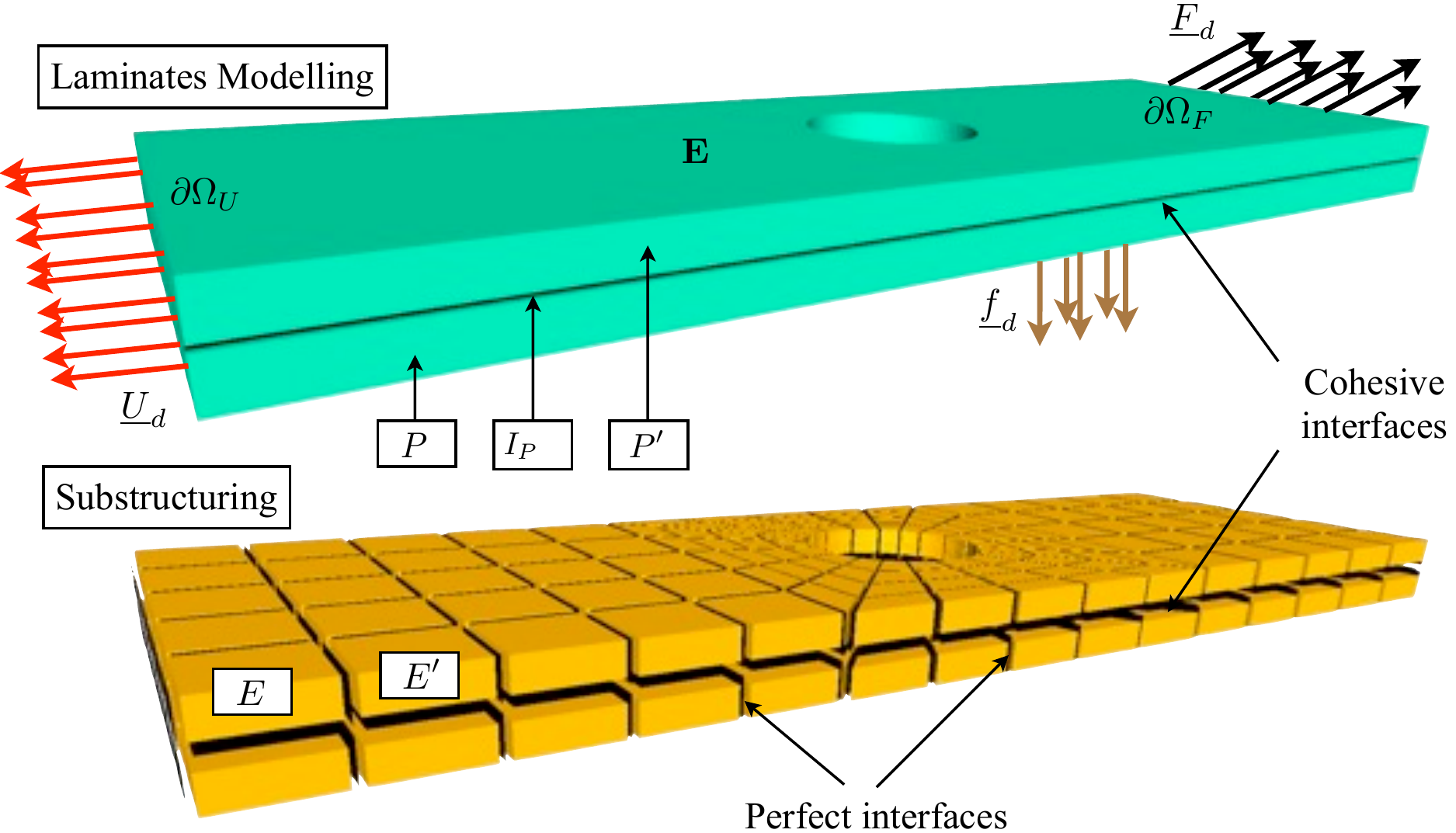}
       \caption{Substructuring of the laminated composite structure}
       \label{fig:decomp_sst_interfaces}
\end{figure}

\begin{figure}
       \centering
       \includegraphics[width=0.65 \linewidth]{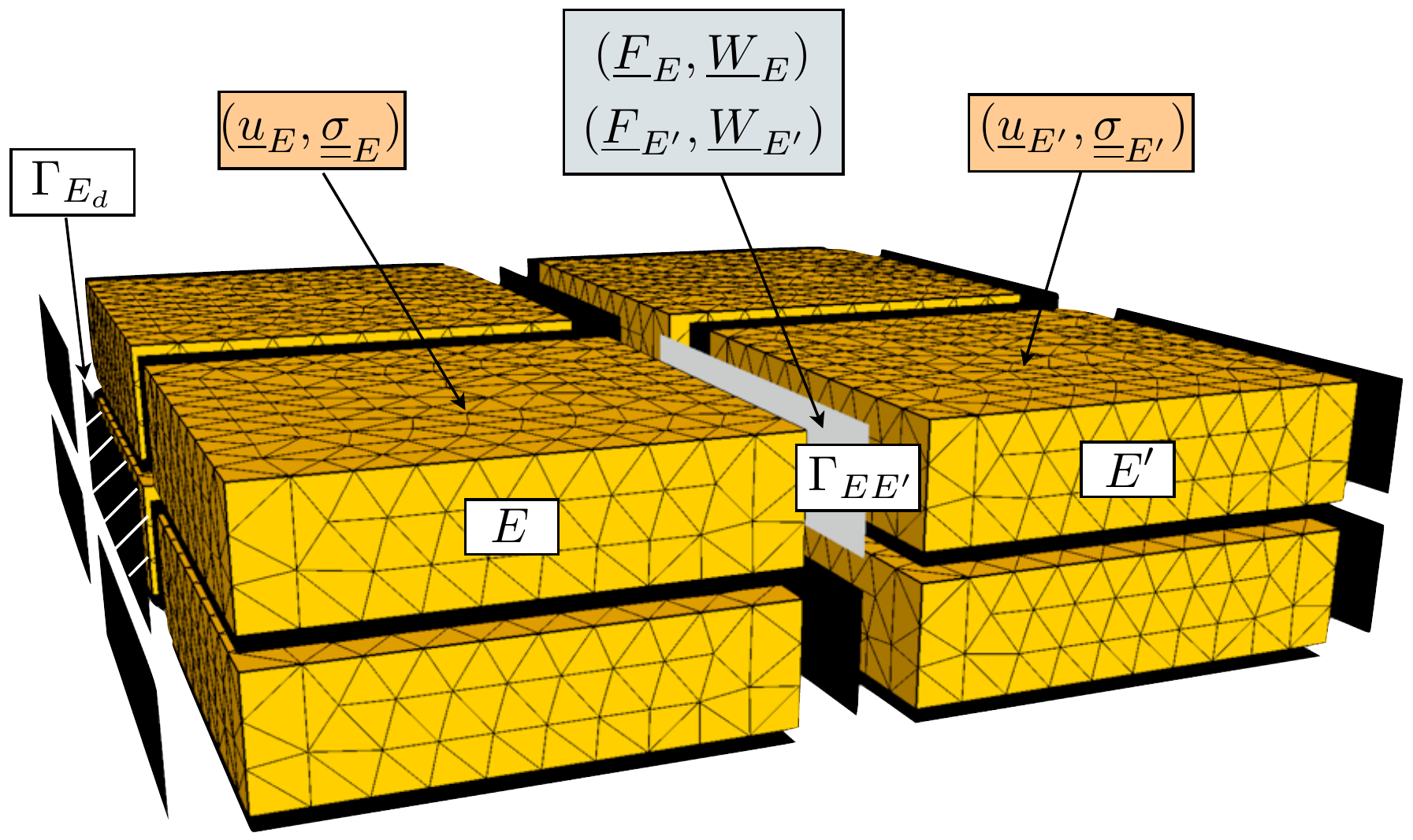}
       \caption{Substructuring of the laminated composite structure}
       \label{fig:champs_interface}
\end{figure}

The laminate structure $\mathbf{E}$ is decomposed into substructures and interfaces as represented in Figure (\ref{fig:decomp_sst_interfaces}). Each of these mechanical entities possesses its own kinematic and static unknown fields linked by its behavior. The substructuring is driven by the will to match domain decomposition interfaces with material cohesive interfaces, so that each substructure belongs to a unique ply  and has a constant linear behavior. Each substructure is defined in a domain $\Omega_E$ such that $E \in \llbracket 1, \ n_E\rrbracket$ ($n_E$ being the total number of substructures) and is connected to a adjacent substructures through interfaces $\Gamma_{EE'}=\partial \Omega_E \cap\partial \Omega_{E'}$ where $E' \in \llbracket 1, \ n_E\rrbracket$ (Figure (\ref{fig:champs_interface})). The surface entity $\Gamma_{EE'}$ applies force distributions $\underline{F}_E$, $\underline{F}_{E'}$ as well as displacement distributions $\underline{W}_E$, $\underline{W}_{E'}$ to Substructure $E$ and Substructure $E'$ respectively. On Substructure $E$ such that $\partial \Omega_{E} \cap \partial \Omega \neq \emptyset$, the boundary condition $(\V{U}_d,\V{F}_d)$ is applied through a boundary interface $\Gamma_{E_d}$. Let us define $\Gamma_E = \bigcup_{E' \in \llbracket 1, \ n_E\rrbracket} \Gamma_{EE'} \cup \Gamma_{E_d}$. We finally introduce $\M{\sigma}_E$, the Cauchy stress tensor, and $\M{\epsilon}(\V{u}_E)$, the symmetric part of the displacement gradient, in substructure $E$.

The substructured quasi-static problem at any computation time $t_{n+1}$ of the time discretization scheme consists in finding $s=(s_E)_{E \in \llbracket 1, \ n_E\rrbracket} $, where $s_E = ( \V{W}_E, \V{F}_E )$, which satisfies the following equations:
\begin{itemize}
\item Kinematic admissibility of Substructure $E$:
\begin{equation}
{\V{u}_E}_{|\Gamma_E} = {\V{W}_{E}}
\end{equation}
\item Static admissibility of Substructure $E$: $\forall ({\underline{u}_E}^\star,{\underline{W}_E}^\star) \in \mathcal{U}_{E} \times \mathcal{W}_{E} \,  / \,  {{\V{u}_E}^\star}_{| \partial \Omega_E} = {{\V{W}_{E}}^\star}, $
\begin{equation}
\label{equation:equilibre_ss}
\displaystyle \int_{\Omega_E} \operatorname{Tr} \left( \M{\sigma}_E \, \M{\epsilon} ({\V{u}_E}^\star) \right) \ d \Omega  = \displaystyle  \int_{\Gamma_E} \V{F}_E . {\V{W}_E}^\star \ d \Gamma  
\end{equation}
\item Linear orthotropic behavior of Substructure $E$:
\begin{equation}
\M{\sigma}_E = K \, \M{\epsilon}(\underline{u}_E)
\end{equation}
\item Behavior of the interfaces $\Gamma_{EE'} \in \Gamma_E$:
\begin{equation}
\label{eq:interface_ref}
\displaystyle \mathcal{R}_{EE'}( \V{W}_{E}  , \V{W}_{E'} , \V{F}_{E} , \V{F}_{E'} ) = 0
\end{equation}
\item Behavior of the interfaces $\Gamma_{{E}_d} \in (\Gamma_E \cap \partial \Omega)$:
\begin{equation}
\label{eq:boundary_ref}
\mathcal{R}_{E_d}( \V{W}_{E}, \V{F}_{E}) = 0 \qquad
(\V{W}_{E}=\V{u}_d \  \text{ on } \partial\Omega_u \ \text{ and } \ \V{F}_{E}=\V{F}_d \ \text{ on }\partial\Omega_f)
\end{equation}
\end{itemize}
 
In delamination analysis, the formal relation $\mathcal{R}_{EE'}=0$ reads:
\begin{itemize}
\item For perfect interface: 
$
\qquad \ \  \, \left\{ \begin{array}{l}
\V{F}_{E} + \V{F}_{E'} = 0 \\
\V{W}_{E} - \V{W}_{E'} = 0
\end{array} \right.
$
\item For cohesive interface:
$
\qquad \left\{ \begin{array}{l}
\displaystyle \V{F}_{E} + \V{F}_{E'} = 0 \\
\displaystyle \V{F}_{E} = \M{K}_{P} \left( \left( \V{W}_{E'} - \V{W}_{E} \right)_{| t \in \llbracket t_0, \  t_{n+1} \rrbracket }  \right) (\V{W}_{E'} - \V{W}_{E})
\end{array} \right.
$ \\where Substructure $E$ (respectively $E'$) belongs to Ply $P$ (respectively $P+1$).
\end{itemize}

\subsection{Iterative resolution of the stable nonlinear substructured problem}
\label{sec:multiscale_resolution}

\begin{figure}[h]
       \centering
       \includegraphics[width=0.45 \linewidth]{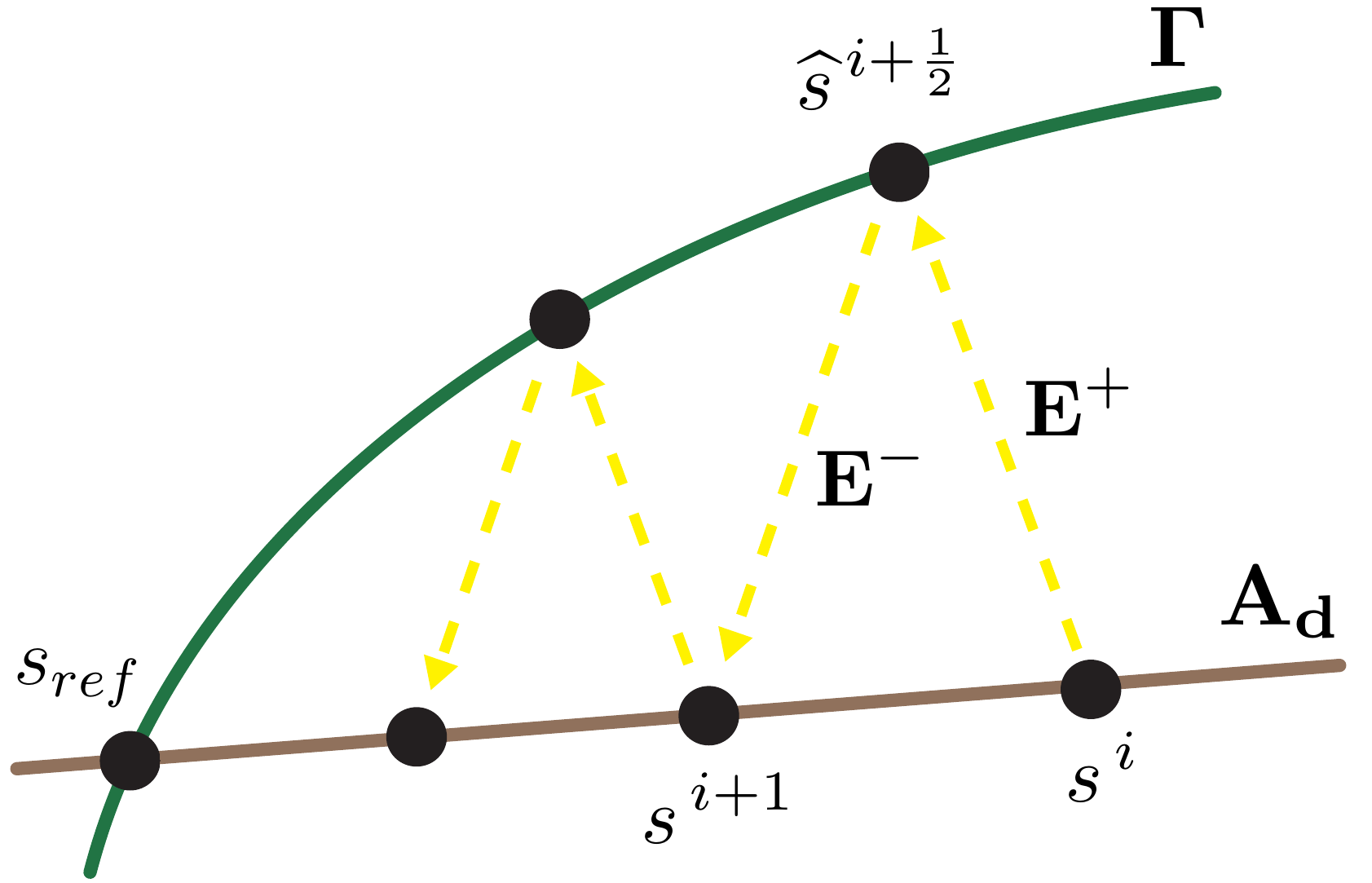}
       \caption{Schematic representation of the LaTIn algorithm}
       \label{fig:latin}
\end{figure}

The equations of the problem can be split into the set of linear equations in substructures (static and kinematic admissibility of the substructures, linear constitutive law of the substructures) and the set of local equations in interface variables (behavior of the interfaces). The solutions $\displaystyle s = (s_E)_{E \in \llbracket 1, \ n_E\rrbracket} = (\V{W}_E ,  \V{F}_E )_{E \in \llbracket 1, \ n_E\rrbracket}$ to the first set of equations belong to Space $\mathbf{A_d}$, while the solutions $\displaystyle \widehat{s} = (\widehat{s}_E)_{E \in \llbracket 1, \ n_E\rrbracket} = (\V{\widehat{W}}_E ,  \V{\widehat{F}}_E )_{E \in \llbracket 1, \ n_E\rrbracket}$ to the second set of equations belong to $\boldsymbol{\Gamma}$. Hence, the converged solution $s_{ref}$ is such that $ s_{ref} \in \mathbf{A_d} \bigcap \boldsymbol{\Gamma} $.

The LaTIn resolution scheme consists in searching for the solution $s_{ref}$ 
alternatively in these two spaces along search directions  $\mathbf{E^+}$ and $\mathbf{E^-}$ (see Fig. \ref{fig:latin}): 
\begin{itemize}
\item 
\textit{Find}  $\widehat{s}^{i+\frac{1}{2}} \in \boldsymbol{\Gamma}  $ \textit{such that}  $ \left( \widehat{s}^{i+\frac{1}{2}} -s^{i} \right) \in \mathbf{E^+} $ (local stage)
\item 
\textit{Find} $s^{i+1} \in \mathbf{A_d}$ \textit{such that}  $\left( s^{i+1} - \widehat{s}^{i+\frac{1}{2}} \right) \in \mathbf{E^-} $ (linear stage)
\end{itemize}
In the following, the subscript $i$ will be dropped.

\paragraph{Local stage}

One searches for a solution $\widehat{s}=(\underline{\widehat{F}}_E,\underline{\widehat{W}}_E)_{E \in \llbracket 1, \ n_E\rrbracket}$ satisfying the local equations on the interfaces ($ \mathcal{R}_{EE'}=0$ or  $\mathcal{R}_{E_d}=0$), and search direction equation $\mathbf{E}^+$, introduced locally on the interfaces :
\begin{equation}
\displaystyle  (\underline{\widehat{F}}_E-\underline{F}_E) - k_E^+ (\underline{\widehat{W}}_E - \underline{W}_E) = 0 
\end{equation}
At this stage, variables $\underline{F}_E$ et $\underline{W}_E$ are known from the previous semi-iteration.

In the case where $ \mathcal{R}_{EE'}=0$ is a nonlinear equation, the local problem is solved by a quasi-Newton algorithm.

\paragraph{Linear stage}

One searches for a solution $s=(\underline{F}_E,\underline{W}_E)_{E \in \llbracket 1, \ n_E\rrbracket}$  verifying the linear equations on each substructure and, at best, a search direction equation $E^-$, local on the interfaces, under the constraint of average equilibrium of the interface forces :
\begin{equation}
\left\{\begin{array}{l}
\displaystyle {\underline{F}_E}_{|\Gamma_{E}} =\displaystyle   \arg\min \left\{ \int_{\Gamma_{E}}\left(\frac{1}{2 \, k_E^-}(\underline{F}_E-\underline{\widehat{F}}_E)^2 +  (\underline{F}_E-\underline{\widehat{F}}_E). (\underline{W}_E - \underline{\widehat{W}}_E)\right) \, d \Gamma \right\} \\
\displaystyle \textrm{under the constraint: } \quad \forall (E',E), \quad \Pi_{|\Gamma_{EE'}}^M \underline{F}_{E_{|\Gamma_{EE'}}}  + \Pi_{|\Gamma_{EE'}}^M \underline{F}_{{E'}_{|\Gamma_{EE'}}} =0
\end{array}\right.
\end{equation}
The macroscopic projectors $\Pi_{|\Gamma_{EE'}}^M$ extract an average of the interface forces, which is transfered into the whole structure. Technically, this stage consists in solving, in parallel, independent linear problems on the sub-structures (using finite elements) and a small macroscopic linear problem which is global over the structure (and discrete by construction).

\subsection{Iterative resolution of the unstable nonlinear problem}

When a snap-back appears in the global behavior of the simulated structure, the incremental LaTin algorithm is switched to a well-known local arc-length algorithm \cite{schellekens93,allix96, kerfriden09}. 
The algebraic nonlinear problem to solve at Time $t_{n+1}$, in an unstable phase, reads: 
\begin{equation}
\label{eq:non-linear-alg}
q_{\textrm{int}} \left(U_{|t_{n+1}},(U_{|\tau})_{\tau < t_{n+1}} \right)  - \lambda_{|t_{n+1}} \, {q_{\textrm{ext}}} = 0
\end{equation} 
The amplitude of the loading $\lambda_{|t_{n+1}}$ is unknown. A control equation inspired from \cite{allix96} is introduced so that the maximum local increment in the jump of displacement over all the cohesive interfaces takes a predefined value $\Delta l$:
\begin{equation}
\label{eq:control}
\mathbf{c}( \Delta U_{|t_{n+1}}) \, \Delta U_{|t_{n+1}} = \Delta l
\end{equation}
where the $\Delta \, . \,$ unknowns are the increments in the quantities over Time step $[t_n \ t_{n+1}]$. $\mathbf{c}$ is then the operator which extract the maximum of the (positive) jump increment.

Classically, the non-linear system (\ref{eq:non-linear-alg}, \ref{eq:control}) is solved by a modified Newton-Raphson scheme: 
\begin{itemize}
\item The linearization of \eqref{eq:non-linear-alg} and \eqref{eq:control} around point $(U^{i},\lambda^{i})$ leads to the system to solve at the prediction step of the $({i+1})^{th}$ iteration of this scheme:
\begin{equation}
\left\{ \begin{array}{l}
\displaystyle \lambda_{|t_{n+1}}^{i+1} = \frac{\Delta l + \mathbf{c}(\Delta U_{|t_{n+1}}^i) \, U_{|t_n}}{\mathbf{c}(\Delta U_{|t_{n+1}}^i) \, \mathbf{K} \left(U_{|t_{n+1}}^i,(U_{|\tau})_{\tau < t_{n+1}} \right)^{-1} {q_{\textrm{ext}}}}\\
\displaystyle U_{|t_{n+1}}^{i+1}  = \lambda_{|t_{n+1}}^{i+1} \, \mathbf{K} \left(U_{|t_{n+1}}^i,(U_{|\tau})_{\tau < t_{n+1}} \right)^{-1} {q_{\textrm{ext}}}
\end{array} \right.
\label{eq:predic_riks}
\end{equation}
The inversion of the linearized stiffness operator (i.e.: the resolution of the linear system $\bar{U} = \mathbf{K} (U_{|t_{n+1}}^i,(U_{|\tau})_{\tau < t_{n+1}} )^{-1} {q_{\textrm{ext}}}$) is performed by using the domain decomposition method described previously (the internal variables of the interfaces are fixed during the resolution)
\item The correction step of the Newton scheme consists in updating Operators $\mathbf{K}$ and $\mathbf{c}$ with respect to the kinematic field $U_{|t_{n+1}}^{i+1}$ found at the prediction stage.
\end{itemize}

\subsection{Stopping criterion}
\label{sec:stop_crit}


\subsubsection{Stable phase (LaTIn algorithm)}
\begin{figure}[htb]
       \centering
       \includegraphics[width=0.45 \linewidth]{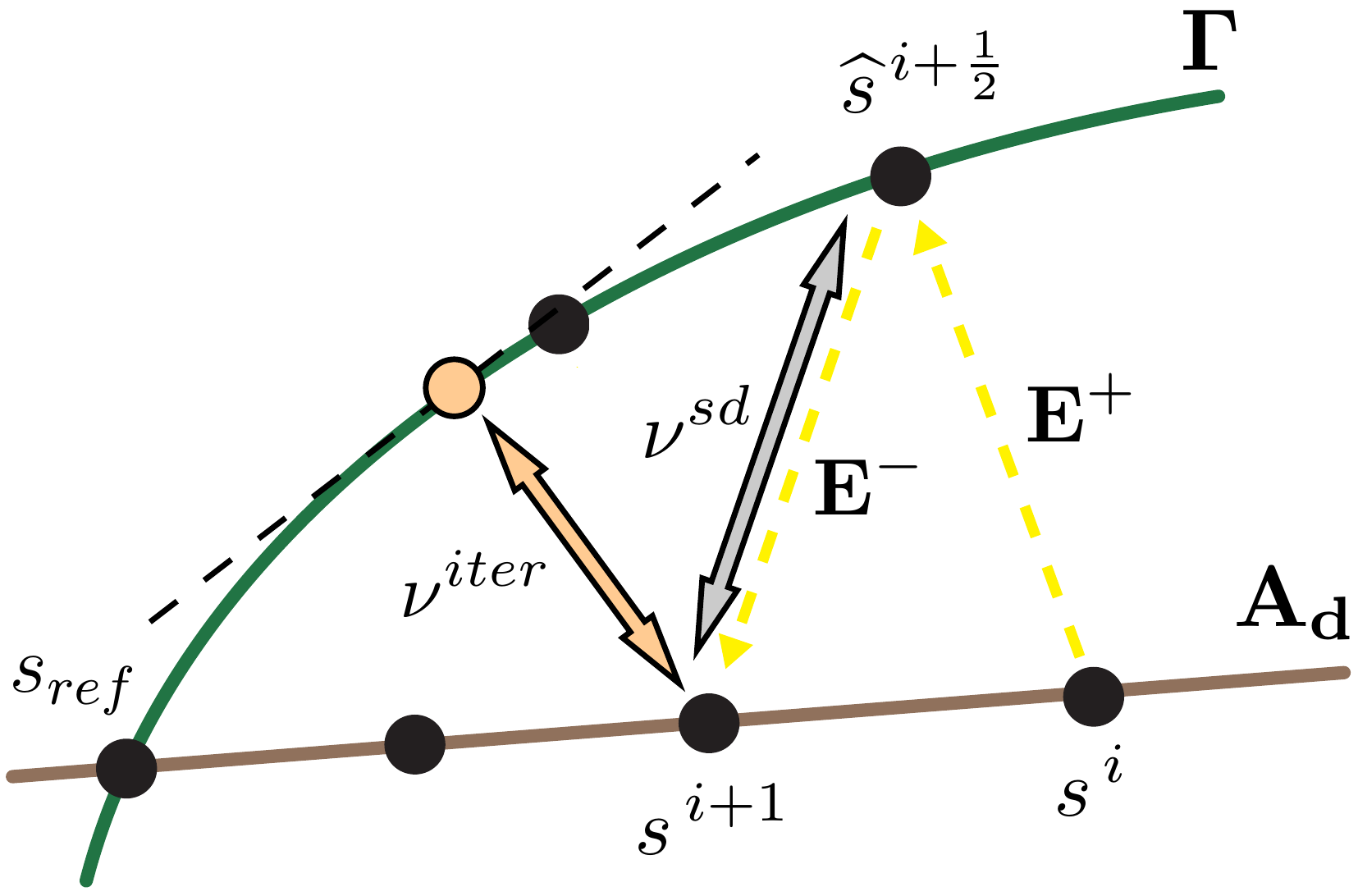}
       \caption{Classical convergence indicator $\nu^{sd}$ of the LaTIn solver and indicator $\nu^{iter}$ based on the error in the constitutive law}
       \label{fig:erreur_rdc2}
\end{figure}
In order to evaluate the convergence of the LaTIn algorithm, one classically measures the distance between spaces $\mathbf{A_d}$ and $\boldsymbol{\Gamma}$ along search direction $\mathbf{E}^-$ \cite{ladeveze99} (criterion labeled $\nu^{sd}$ on Figure (\ref{fig:erreur_rdc2}), ``\it sd\normalfont '' standing for ``search direction''). In the work of \cite{kerfriden09,ladeveze09}, a new criterion based on the error in the constitutive law has been successfully assessed (in order, originally, to get rid of the dependency of Convergence indicator $\nu^{sd}$ on the parameters of the LaTIn solver). The solutions resulting from a linear stage of the LaTIn algorithm satisfy all the equations of the substructured reference problem except the interfaces laws \eqref{eq:interface_ref} and \eqref{eq:boundary_ref}, whose residuals can be easily computed (Figure (\ref{fig:erreur_rdc2})). More precisely, a solution $s^{i+1} \in \mathbf{A_d}$ being reached, an indicator of the convergence of the algorithm is given by integrating the corresponding local residuals of the interface behaviors over the structure (residuals of Equations \eqref{eq:interface_ref} and \eqref{eq:boundary_ref} evaluated for $\displaystyle s^{i+1} = (s^{i+1}_E)_{E \in \llbracket 1, \ n_E\rrbracket} = (\V{W}^{i+1}_E ,  \V{F}^{i+1}_E )_{E \in \llbracket 1, \ n_E\rrbracket}$). 

Let us call $q$ the number of interface relations being used (i.e.: the number of distinct interface behaviors $(\mathcal{R}_{EE'}=0)_{(E,E') \in \llbracket 1, \ n_E\rrbracket^2}$ and $(\mathcal{R}_{E_d}=0)_{E \in \llbracket 1, \ n_E\rrbracket}$ used in the structure). In our case, $q=4$ (perfect interfaces, cohesive interfaces with homogeneous constants, Dirichlet and Neumann boundary conditions). ${\bar{\Gamma}}_i$ is the set interfaces of Behavior $i$, for all $i \in  \llbracket 1, \ q \rrbracket $. The vectorial relations $\mathcal{R}_{EE'}=0$ for $i \in  \llbracket 1, \ q \rrbracket $ and ${\Gamma_{EE'} \in \bar{\Gamma}_i \ \text{or} \ \Gamma_{E_d} \in \bar{\Gamma}}_i$ are made of $p_i$ vectorial equations $\mathcal{Q}_{ij}=0$ (2 equations for cohesive or perfect interfaces in 3D, 1 equation for boundary interfaces). Here, subscript j ranges from $1$ to $p$.  Convergence indicator $\nu^{iter}$ (``\textit{iter}'' stands for ``iterative'') reads:
\begin{equation}
\label{eq:nuiter}
\displaystyle \left( {\nu^{iter}} \right)^2 = \sum_{i=1}^q \sum_{j=1}^{p_i} \left( {\nu^{iter}_{ij}} \right)^2 
 \qquad   \textrm{where} \qquad 
\left( {\nu^{iter}_{ij}} \right)^2 = \frac{ \displaystyle \sum_{\Gamma \in \bar{\Gamma}_i}   \int_{\Gamma} {\mathcal{Q}_{ij}} .  \mathcal{Q}_{ij} \ d \Gamma }{  \displaystyle \sum_{\Gamma \in \bar{\Gamma}_i}   \int_{\Gamma} \mathcal{\widetilde{Q}}_{ij} .  \mathcal{\widetilde{Q}}_{ij} \ d \Gamma }
\end{equation}
where, in the case of delamination (i.e : involving perfect and cohesive LaTIn interfaces):
\begin{itemize}
\item on a perfect interface $\GammaEEp \in \bar{\Gamma}_1$ :
\begin{equation}
\begin{array}{l}
\displaystyle \mathcal{Q}_{11} = \V{F}_E + \V{F}_{E'} \\
\displaystyle \mathcal{\widetilde{Q}}_{11} = \V{F}_E - \V{F}_{E'} 
\end{array}
 \qquad \qquad
\begin{array}{l}
\displaystyle \mathcal{Q}_{12} = \V{W}_E - \V{W}_{E'} \\
\displaystyle \mathcal{\widetilde{Q}}_{12} = \V{W}_E + \V{W}_{E'}
\end{array}
\end{equation}
\item on a cohesive interface $\GammaEEp \in \bar{\Gamma}_2$ :
\begin{equation}
\label{eq:cohesive_time_independent}
\begin{array}{l}
\displaystyle \mathcal{Q}_{21} = \V{F}_E - \M{K}_{P}\left((\V{W}_{E'}-\V{W}_{E})_{|t \in \{ t_{n+1} , [ 0 \ t_{n} ] \} }\right) \left( \V{W}_{E'} - \V{W}_{E} \right) \\
\displaystyle \mathcal{\widetilde{Q}}_{21} = \V{F}_E + \M{K}_{P}\left((\V{W}_{E'}-\V{W}_{E})_{|t \in \{ t_{n+1} , [ 0 \ t_{n} ] \} } \right) \left( \V{W}_{E'} - \V{W}_{E} \right) \\
\displaystyle \mathcal{Q}_{22} = \V{F}_{E'} - \M{K}_{P}\left((\V{W}_{E'}-\V{W}_{E})_{|t \in \{ t_{n+1} , [ 0 \ t_{n} ] \} }\right) \left( \V{W}_{E}-\V{W}_{E'} \right)  \\
\displaystyle \mathcal{\widetilde{Q}}_{22} = \V{F}_{E'} + \M{K}_{P}\left((\V{W}_{E'}-\V{W}_{E})_{|t \in \{ t_{n+1} , [ 0 \ t_{n} ] \} } \right) \left( \V{W}_{E}-\V{W}_{E'} \right)
\end{array}
\end{equation}
where $P \in \llbracket 1, \ N_P-1\rrbracket$. Note that, in Equation \eqref{eq:cohesive_time_independent}, the history of the interface variables during the current load increment is not taken into account, for it is unknown at this stage of the resolution procedure.
\item on an interface transmitting Neumann's boundary condition $\GammaEd \in \bar{\Gamma}_3$ :
\begin{equation}
\displaystyle \mathcal{Q}_{31} = \V{F}_E - \V{F}_{d} \qquad \qquad
\displaystyle \mathcal{\widetilde{Q}}_{31} = \V{F}_E+ \V{F}_{d} 
\end{equation}
\item on an interface transmitting Dirichlet's boundary condition $\GammaEd \in \bar{\Gamma}_4$ :
\begin{equation}
\displaystyle \mathcal{Q}_{41} = \V{W}_E - \V{W}_{d}  \qquad \qquad
\displaystyle \mathcal{\widetilde{Q}}_{41} = \V{W}_E + \V{W}_d
\end{equation}
\end{itemize}
The computation of this criterion is very cheap as it simply requires local integration over each interface of the domain decomposition method, and a global sum of these local contributions over the structure

\subsubsection{Unstable phase (arc-length procedure)}

The convergence of the algorithm is evaluated by computing Criterion $\nu^{iter}$ after each prediction stage of the Newton scheme (the residual of the control equation, which has no physical meaning, is not accounted for).

\subsubsection{Typical values}

Our experiments of delamination analysis within the LaTIn framework have shown that a stopping criterion $\nu^{iter}$ set to ${\nu^{iter}_d} = 1 \times 10^{-2}$ (``\textit{d}'' stands here for ``desired'') is usually sufficient to ensure a global convergence of the iterative process (at least, crack fronts are correctly localized, which means that the large wavelength piece of information is correctly captured). In our simulations, and in order to force an accurate convergence of the local quantities (equilibrium of the interface forces and verification of the cohesive law in the process zones), ${\nu^{iter}_d}$ is set to $1 \times 10^{-3}$. \\

%% file: chap_3.tex

In this section, we combine the ideas detailed in Section \ref{sec:criterion} to estimate the time discretization error, and the developments of the last section, dedicated to the evaluation of the convergence of the iterative parallel resolutions to derive a new time discretization error indicator, suited (but not restricted) to the mixed domain decomposition strategy. Based on this new indicator, an automatic procedure to control the load increments is  derived.

\subsection{Time discretization error criterion in a domain decomposition framework}

\subsubsection{Definition}

A sufficiently converged solution of the reference problem being reached at Current time $t_{n+1}$, by making use of the LaTIn-based resolution strategy, a continuous solution is constructed over $[ t_n \ t_{n+1} ]$, as described in Section \ref{sec:criterion}. A new time discretization error criterion ${\nu^{time,dd}}$ (``\textit{dd}'' stands for ``domain decomposition'') is computed with respect to the interpolated solution:
\begin{equation}
\nu^{time,dd}_{|t_{n+1}} = \max_{i \in  \llbracket  0, \ N_s  \rrbracket} \nu^{interp,dd}_{| \bar  t_i}
\end{equation}
where we recall that $N_s+1$ is the number of intermediate times $(\bar t_i)_{i \in  \llbracket  0, \ N_s  \rrbracket}$ such that $t_n \leq \bar t_i \leq t_{n+1}$ at which intermediate solutions are constructed, and Criterion $\nu^{interp,dd}$ is evaluated. 

$\nu^{interp,dd}_{| \bar  t_i}$ is computed by the same formulas defining $\nu^{iter}_{| \bar  t_i}$, except that the history of the interface variables is known ``continuously'' over Time interval $[ 0 \  \bar  t_i ]$ (from the interpolation):
\begin{equation}
\label{eq:nuiter1}
\displaystyle \left( {\nu_{| \bar  t_i}^{interp,dd}} \right)^2 = \sum_{i=1}^q \sum_{j=1}^p \left( {\nu^{interp,dd}_{ij_{| \bar  t_i}}} \right)^2 
 \quad   \textrm{where} \quad 
\left( {\nu^{interp,dd}_{ij}} \right)^2 = \frac{ \displaystyle \sum_{\Gamma \in \bar{\Gamma}_i}   \int_{\Gamma} {\mathcal{Q}_{ij_{| \bar  t_i}}} .  \mathcal{Q}_{ij} \ d \Gamma }{  \displaystyle \sum_{\Gamma \in \bar{\Gamma}_i}   \int_{\Gamma} {{\mathcal{\widetilde{Q}}_{ij}}}.  \mathcal{\widetilde{Q}}_{ij} \ d \Gamma }
\end{equation}
and, in Equation \eqref{eq:cohesive_time_independent}, the stiffness operator of the cohesive interfaces ${\Gamma_{EE'}}_{(E,E')\in \llbracket 1, \ n_E\rrbracket^2}$ is replaced by the reconstructed operator $\M{K}_{P}\left((\V{W}_{E'}-\V{W}_{E})_{|t \in [0 \ \bar  t_i] } \right)$.

Note that Time discretization error criteria ${\nu^{time,dd}}$ and ${\nu^{time}}$ are slightly different. In Section \ref{sec:criterion}, we assumed that the solutions obtained at Times $t_n$ and $t_{n+1}$ satisfied the global static admissibility (Equation \eqref{eq:static_admissibility}). This assumption cannot be made anymore if the solver used is the mixed domain decomposition strategy (unless the convergence criterion threshold is set to a very low value, which would be ineffective, from a numerical point of view). Indeed, the equilibrium is only satisfied in terms of substructures and macroscopic interfaces variables. In addition, the kinematic admissibility is not fully satisfied, for each ply has been decomposed into substructures which are imperfectly bonded during the resolution process. Hence, the new time discretization error criterion ${\nu^{time,dd}}$ measures not only the cohesive law residual, but also an interface microscopic equilibrium residual (which is small) and a jump of displacement through perfect interfaces, both due to a partial convergence of the iterative solver. 

\begin{figure}[htb]
       \centering
       \includegraphics[width=0.75 \linewidth]{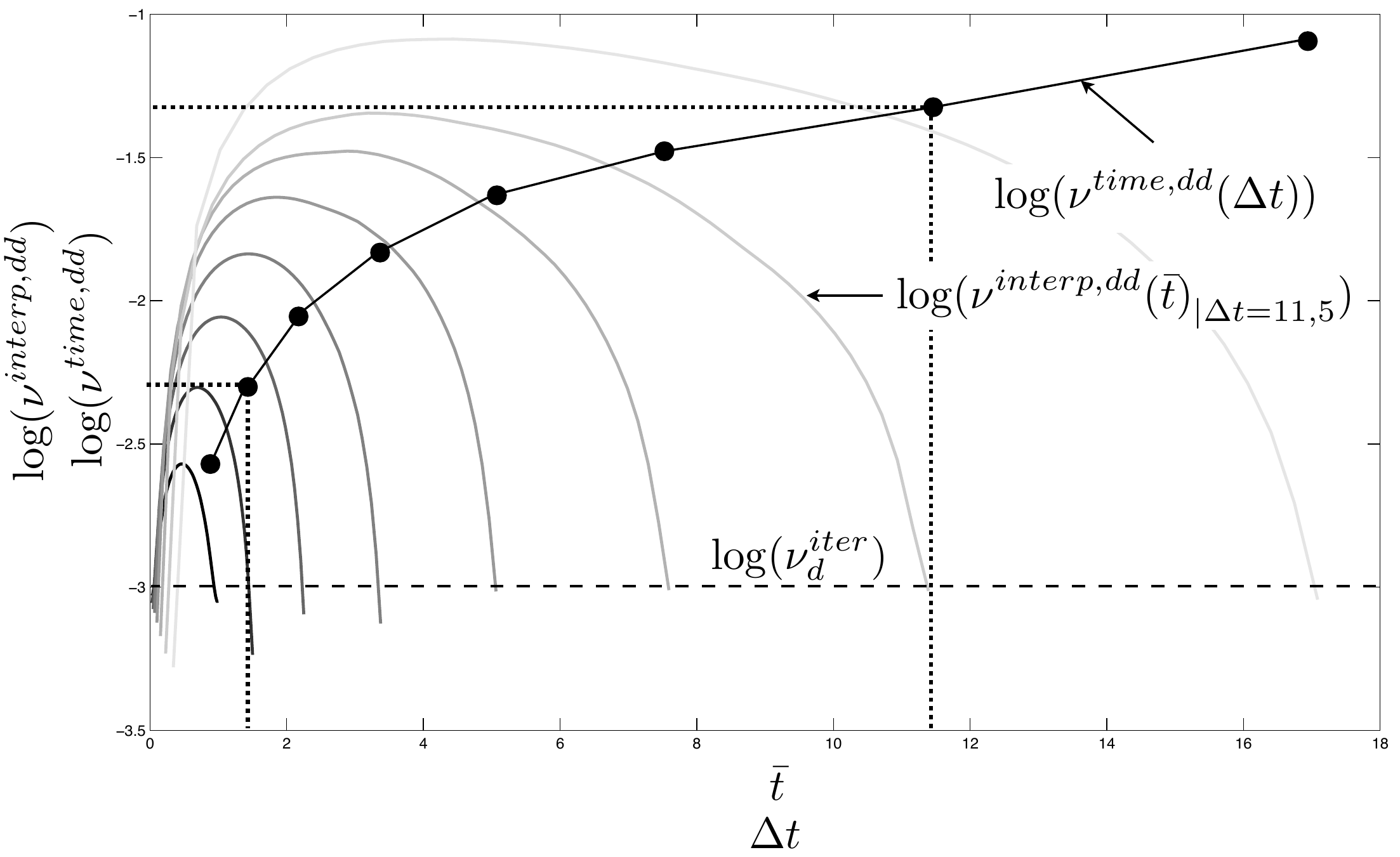}
       \caption{Grey curves : evolution of $\nu^{interp,dd}$ as a function of $\bar t \in [t_n \ t_{n+1}]$ for different values of $\Delta t$. Black curve: Evolution of $\nu^{time,dd}$ with respect to $\Delta t$ (maximum values of the gray curves)}
       \label{fig:different_dl}
\end{figure}

\subsubsection{Properties} 
\label{sec:toto}

Figure (\ref{fig:different_dl}) shows the evolution of $\nu^{interp,dd}$ within a time interval $ [t_n \ t_{n+1}] $ for a given computation time $t_n$ of the unstable delamination simulation represented Figure (\ref{fig:plaque_definition}) (which will be detailed later on). 
The different gray curves correspond to different values of the time increment $\Delta t$ (value of the prescribed arc-length in this case). 
Note that the the value of ${\nu^{interp,dd}}$ at Computation time $t_n$ is the value $\nu^{iter}_d$ of $\nu^{iter}$ which has been prescribed to ensure a sufficient convergence of the LaTIn iterative process. 
From this set of parametric studies, the values of $\nu^{time,dd}$ can be plotted with respect to $\Delta t$ (maximum values of $\nu^{interp,dd}$ over $ [t_n \ t_{n+1}] $, black points on the figure). The resulting function (black interpolated curve) is monotonic. 

One can also remark that even when a large time step is prescribed, the curve $\nu^{interp,dd}$ as a function of $\bar t \in [t_n \ t_{n+1}]$ is smooth. Thus, a relatively small number of evaluation of this residual over the time step is sufficient to obtain an accurate value of the time discretization error criterion $\nu^{time,dd}$. In addition, as the computation of Criterion $\nu^{interp,dd}$ is cheap, even a large number of intermediate time steps would not alter the numerical efficiency of the strategy. In practice, we choose  $N_s=10$.

\subsection{Adaptive load increment procedure}


Our aim is to solve the delamination problem at Computation time $t_{n+1}$ under the constraint of given value $\nu^{time,dd}_d$ of the time discretization error indicator, the current time increment $\Delta t=t_{n+1}-t_n$ (\textit{i.e.}: the prescribed arc-length or the load increment) being set as a new unknown.
A quasi-Newton technique is used to solve the nonlinear constraint equation:
\begin{equation}
\label{eq:prescribed_nu}
\nu^{time,dd}(\Delta t) - \nu^{time,dd}_d = 0
\end{equation}
Basically, each iteration of this scheme is decomposed in three steps:
\begin{itemize}
\item a linear step, where a value of the time increment is predicted (see formulas (\ref{eq:prediction},\ref{eq:prediction2}) in next paragraph).
\item a correction stage where the full delamination problem is solved, at the current time step $t_{n+1}$, until convergence of the underlying nonlinear solver. The time increment $\Delta t$ is here fixed to its predicted value.
\item the computation of a convergence indicator (norm of the residual of Equation \eqref{eq:prescribed_nu})
\end{itemize}

\begin{figure}[ht]
       \centering
       \includegraphics[width=0.65 \linewidth]{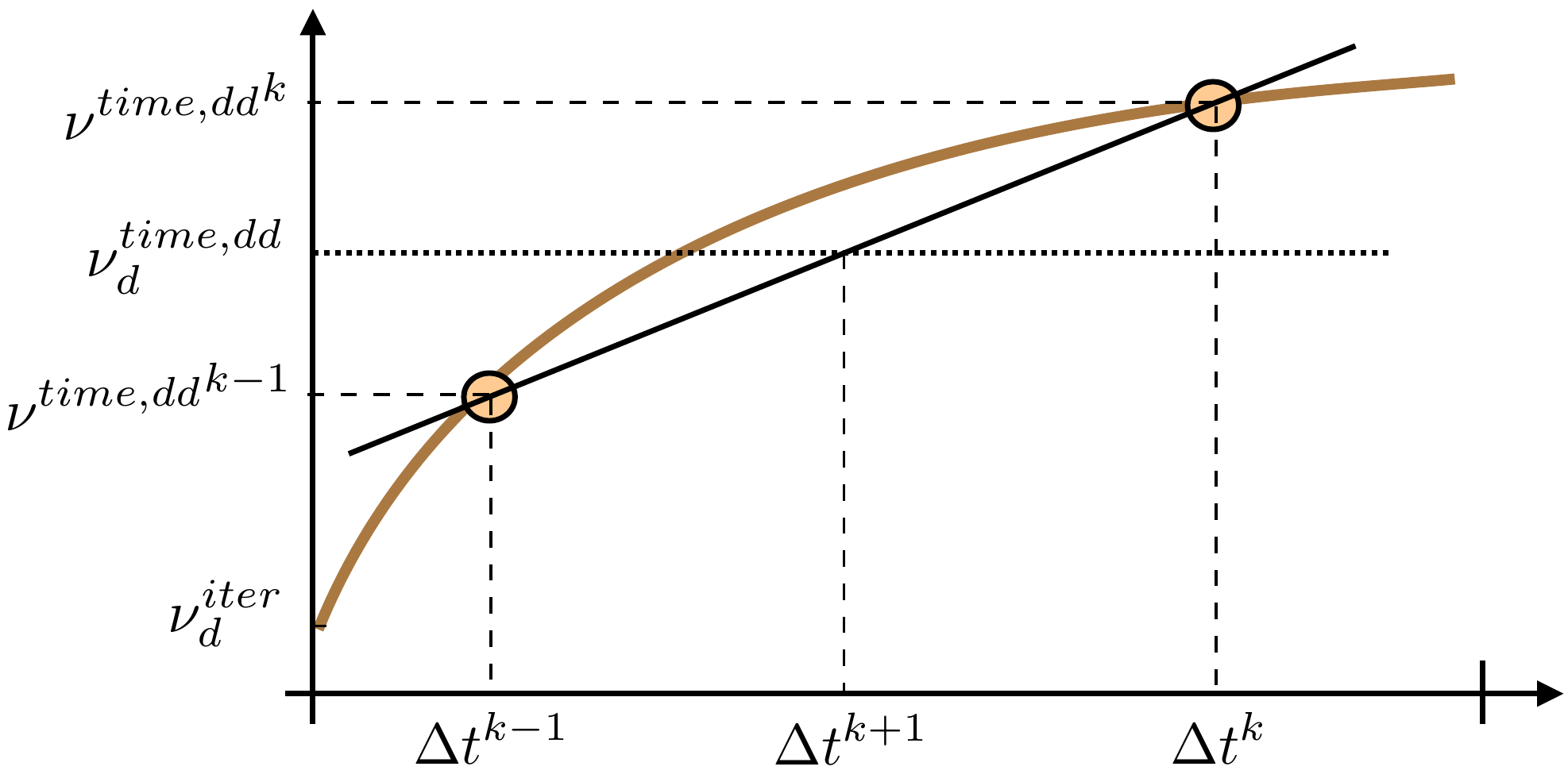}
       \caption{Prediction step of the Newton algorithm designed to solve the delamination problem under the constraint of given time discretization error criterion}
       \label{fig:implicite_control}
\end{figure}

The linear prediction stage at Iteration $k+1$ of Computation time $t_{n+1}$ consists in solving the linearized relation linking $\nu^{time,dd}$ to the time increment $\Delta t$ (see Figure (\ref{fig:implicite_control})). This prediction is done by the following formula: 
\begin{equation}
\label{eq:prediction}
\displaystyle \Delta t^{k+1} = \Delta t^{k-1} +  \frac{ \nu^{time,dd}_d - {\nu^{time,dd}}^{k-1}  }{  {\nu^{time,dd}}^k- {\nu^{time,dd}}^{k-1}  } \left( \Delta t^{k} - \Delta t^{k-1} \right)
\end{equation}
Previous formula does not warranty the prediction of a positive arc-length (the function to linearize is concave). If a negative time increment is predicted, equation \eqref{eq:prediction} is replaced by the following relation, which has empirically shown good convergence properties:
\begin{equation}\label{eq:prediction2}
\Delta t^{k+1} = \sqrt{ \frac{ \nu^{time,dd}_d }{ {\nu^{time,dd}}^k } } \ \Delta t^{k} 
\end{equation}

%% file: chap_4.tex
\subsection{First numerical studies of the time step control procedure on the stable DCB case}

Let us fully detail the results obtained on the stable DCB case presented in Section \ref{sec:reference_problem}, Figure \ref{fig:DCB}. This problem is globally stable and solved using, at each computation time, the LaTIn-based domain decomposition strategy. In the four simulations presented Figure \ref{fig:DCB_time}, the time increments have been obtained by prescribing increasing values of $\nu_{d}^{time,dd}$, respectively $1.0 \times 10^{-2}$, $5.0 \times 10^{-2}$,  $2.0 \times 10 ^{-1}$ and $3.5 \times 10^{-1}$. The resulting average time increment increases, the total number of computation times $N$ being respectively equal to 69, 21, 9 and 5. 

As explained in Section \ref{fig:DCB_time}, the damage state in the cohesive interfaces tends to the one obtained for very small load increments (case 1) when the value of $\nu_{d}^{time,dd}$ decreases. More precisely, the delaminated area of the cohesive interfaces (\textit{i.e.}: the dissipated energy) converges in a monotonic manner with decreasing values of threshold of the time discretization indicator. When this threshold is set to a value smaller than  $2.0 \times 10 ^{-1}$, the error made in terms of dissipated energy is not significant.

Though, this test case is too specific (stable, only one crack front) to give a reliable threshold value $\nu_{d}^{time,dd}$ which should be applied in the general case in order to insure a sufficient convergence of the solution with respect to time.


\subsection{Unstable holed-plate delamination problem}

\begin{figure}[htb]
       \centering
       \includegraphics[width=1. \linewidth]{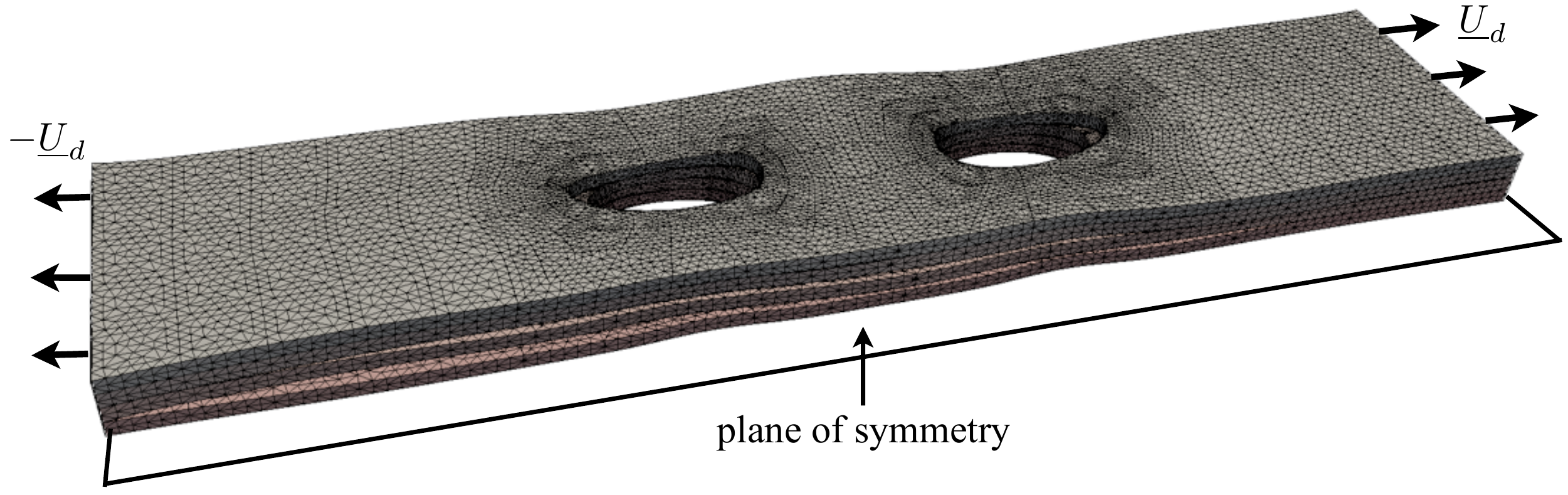}
       \caption{Definition of the holed plate problem (317 000 d.o.f.)}
       \label{fig:plaque_definition}
\end{figure}

We consider a eight-plies holed-plate structure, under traction (prescribed displacements). The plies are orthotropic (stiffness ratio: $1/20$) and the stacking sequence is $[0 \, \pm 45 \ 90]_s$, which leads to the initiation of delamination due to edge effects. The initial stiffness properties of the cohesive interfaces are obtained by the same homogenization in the ``thickness'' of the interfaces (one tenth of the thickness of the plies) which has been described for the DCB problem in Section \ref{sec:reference_problem}. Due to the material and structural symmetries, only the top half of the structure is simulated. The unstable quasi-static time analysis is performed by making use of the arc-length procedure described in Section \ref{sec:latin}. 
The global response curve (plotted in Figure (\ref{fig:2_holes_time_comp}), Case 3) of this case shows two main zones of instability. The first one corresponds to a an unstable propagation of the delamination in the $[-45 \ +45]$ interface while the second one is a crack propagation in the $[0 \ - 45]$ interface, both mainly in shear mode.

\subsubsection{Prescribed time step (Cases numbered 1, 2 and 3 in the whole analysis of the results)}
\begin{figure}[htb]
       \centering
       \includegraphics[width=1. \linewidth]{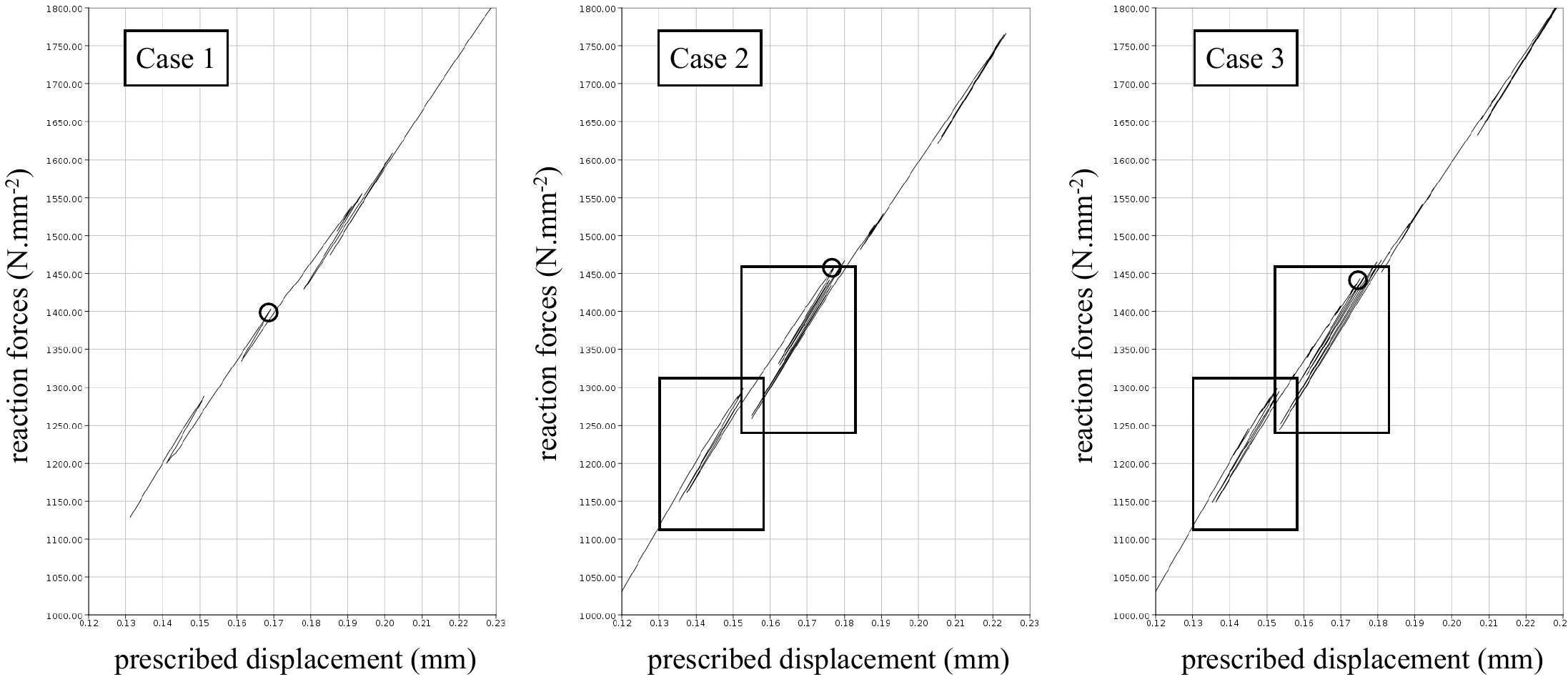}
       \caption{Global reaction force versus prescribed displacement in the holed plate problem under traction, three different predefined arc-lengths being applied (Cases 1,2 and 3)}
       \label{fig:2_holes_time_comp}
\end{figure}

\begin{figure}[htb]
       \centering
       \includegraphics[width=1. \linewidth]{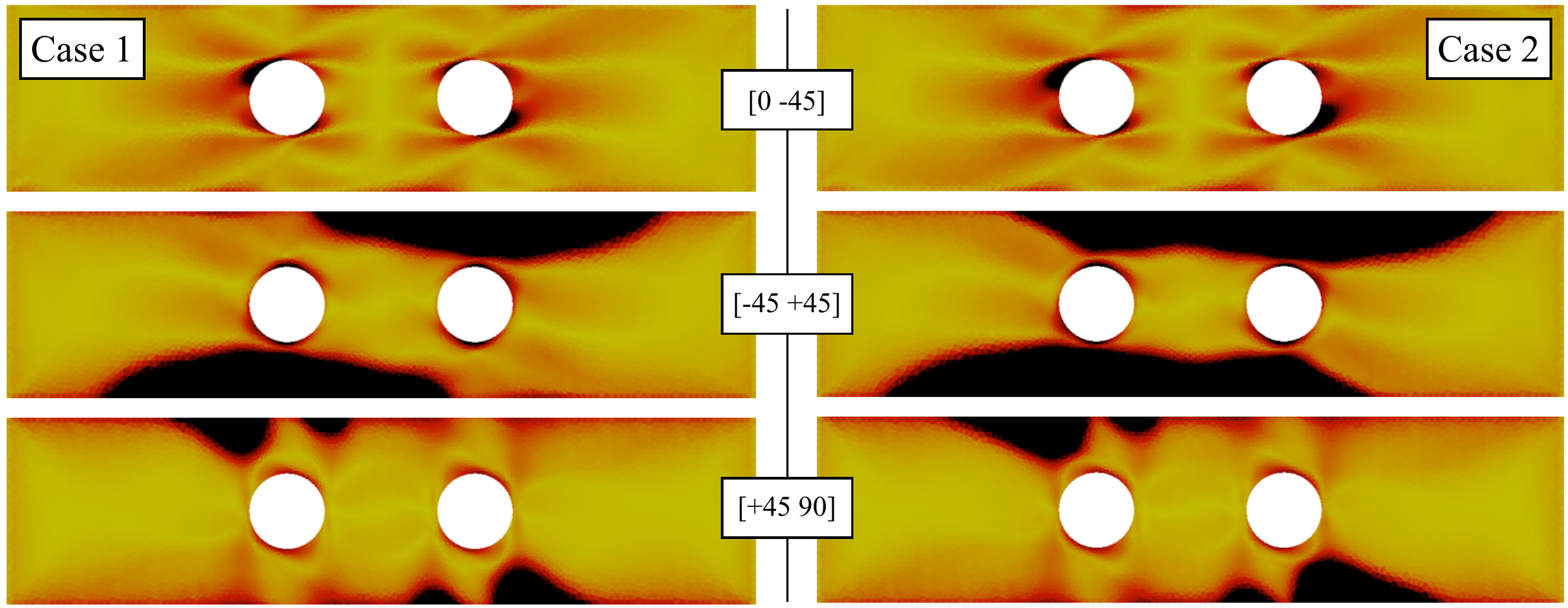}
       \caption{Damage state in the interfaces of the holed plate at the beginning of a global instability in the case of a coarse time grid (Case 1), and in a converged case (Case 2). A Fixed arc-length is prescribed in both cases. In the first case, the damage in the $[-45 \,+45]$ interface is underestimated. }
       \label{fig:2_holes_time}
\end{figure}

The first set of simulations is performed by successively prescribing three different fixed arc-lengths. The arc-length which has been arbitrarily chosen in Case 1 is divided by three in Case 2, and by nine in Case 3. Instabilities appear in the global response of the structure (Figure (\ref{fig:2_holes_time_comp})). 
Figure (\ref{fig:2_holes_time}) shows the damage maps in the $[0 \ - 45]$ $[-45 \ +45]$ and $[-45 \ 90]$ cohesive interfaces in a monotonic phase of the global behavior (limit point after which the delamination front evolves in an unstable manner in the $[0 \, - 45]$ interface, which corresponds to the circled point in all graphs of Figure (\ref{fig:2_holes_time_comp})). 
This particular point of interest has been reached respectively in 8, 63 and 237 time increments. Note that we do not aim at discussing the validity of the solutions reached but at ensuring that the incremental strategy follows the equilibrium path of the converged solution with respect to the time.

\begin{figure}[htb]
       \centering
       \includegraphics[width=0.65 \linewidth]{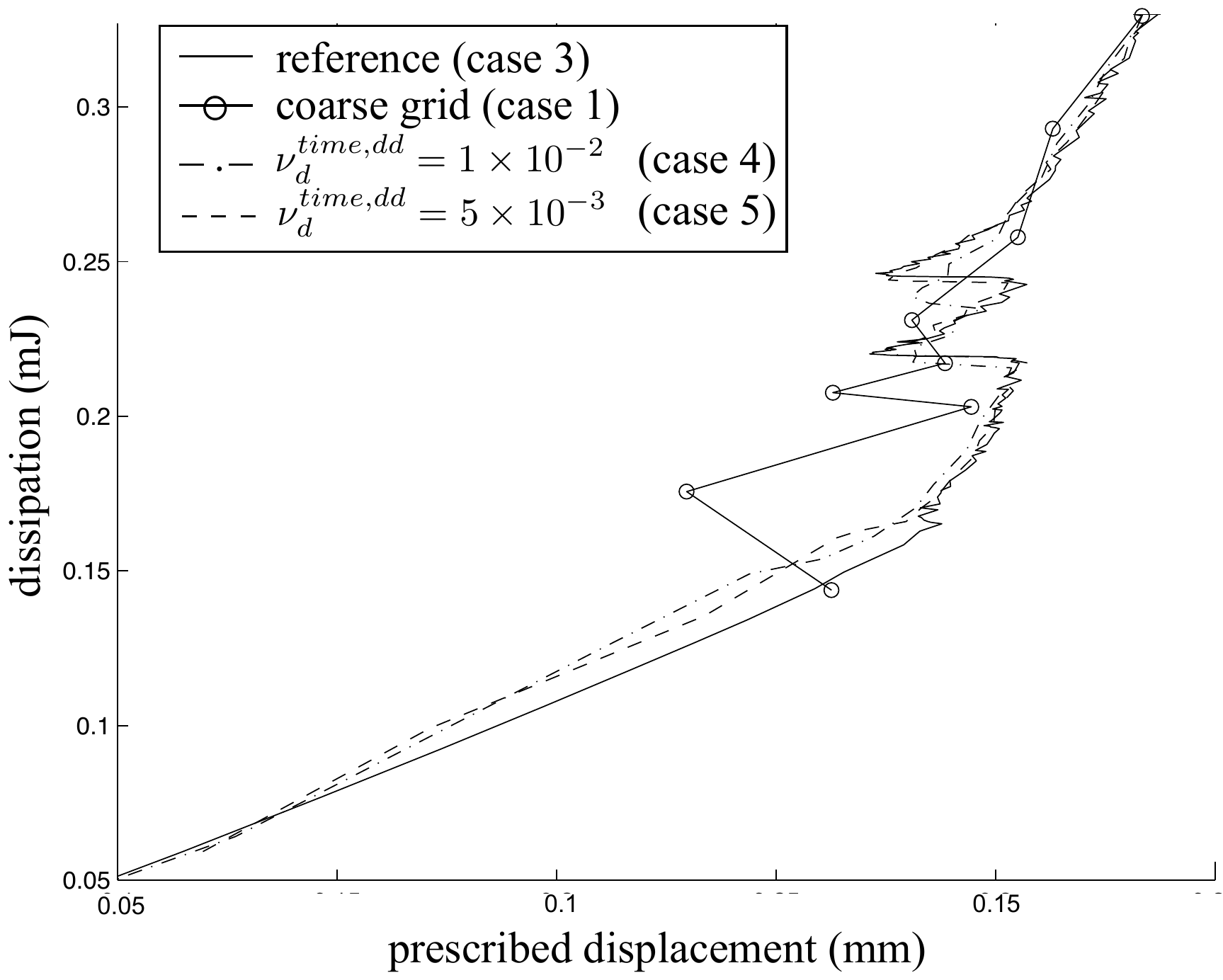}
       \caption{Dissipation versus loading curves for different resolution strategies: explicit fine and coarse time steps (Cases 1 and 3, fixed arc-length) or automatically controlled time increments (Cases 4 and 5, fixed time discretization error)}
       \label{fig:dissi}
\end{figure}

No significant difference can be observed in the damage maps and global response curves corresponding to the two finest analysis in time, which means that the solutions are sufficiently converged with respect to time in Cases 2 and 3. In Case 1, the time increments are too coarse, which results in the incremental resolution procedure to follow a different equilibrium path (see the damage maps in Figure (\ref{fig:2_holes_time})). This phenomenon can impair the global response of the structure, as it can be seen on Figure (\ref{fig:2_holes_time_comp}). The instability phases framed on the converged solutions (Cases 2 and 3) are wrongly predicted in Case 1. These differences appear even more clearly on the dissipated energy versus prescribed displacement curves plotted on Figure (\ref{fig:dissi}) (the curves labeled ``reference'' and ``coarse grid'' refer respectively to Cases 3 and 1), corresponding to the first global instability and to the following stable phases of the time analysis. 
\begin{figure}[htb]
       \centering
       \includegraphics[width=0.75 \linewidth]{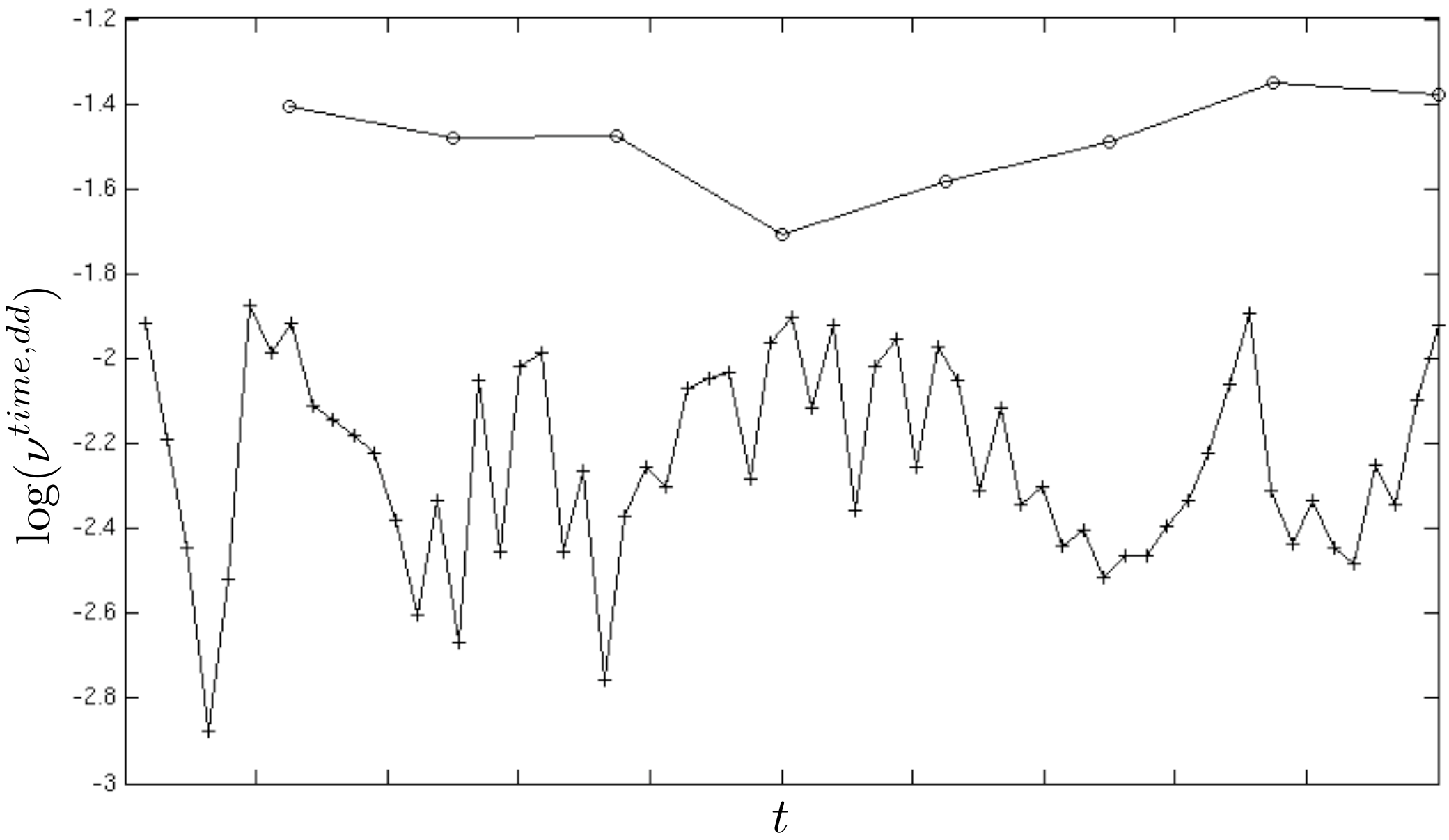}
       \caption{Time discretization error criterion as a function of the Computation times in Case 1 (circles, coarse time steps)  and Case 2 (crosses, small time steps). A Fixed arc-length is prescribed in both cases.}
       \label{fig:2_holes_time_err}
\end{figure}

Figure (\ref{fig:2_holes_time_err}) presents the values of $(\nu^{time,dd}_{|t_n} )_{n \in \llbracket 1, \ \bar{N}\rrbracket}$ as a function of the successive computation times in Cases 1 and 2, from the beginning of the analysis to the starting point of the second global instability. One can see that in Case 2, in which the time increments are sufficiently small to let the iterative algorithm follow the correct equilibrium path, the values of the discretization error indicator $\nu^{time,dd}$ range from $1 \times 10^{-3}$ to $1 \times 10^{-2}$. Conversely, we show in the next set of studies that setting the threshold value $\nu^{time,dd}_d$ of the time control procedure to the maximum of the values $\nu^{time,dd}$ obtained in Case 2  permits to obtain a correctly predicted solution.

\subsubsection{Control of the time step (Cases numbered 4 and 5)}

\begin{figure}[htb]
       \centering
       \includegraphics[width=1. \linewidth]{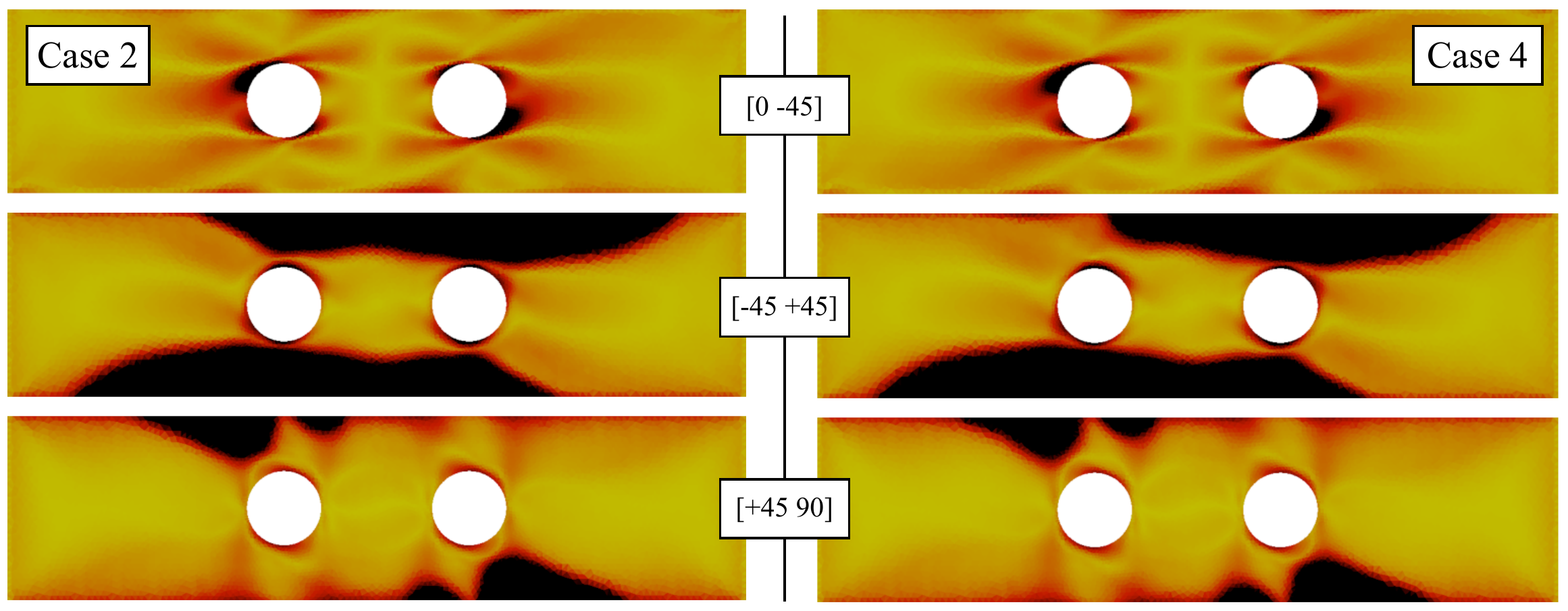}
       \caption{Damage state in the interfaces of the holed plate at the beginning of an instability obtained in a converged case (prescribed arc-length, Case 2), and by using the time control procedure (fixed discretization error, Case 4)}
       \label{fig:2_holes_time_control}
\end{figure}

The second set of simulations  is performed by making use of the procedure described in Section \ref{sec:time_control} to control the successive prescribed arc-length. $\nu_{d}^{time,dd}$ is successively set to  $1 \times 10^{-2}$ (Case 4) and $5 \times 10^{-3}$ (Case 5). The damage state in the cohesive interfaces at the beginning of the second global instability phase predicted by prescribing  $\nu_{d}^{time,dd} = 1 \times 10^{-2}$ is very closed to the one obtained in Case 2 of our first set of simulations (see Figure (\ref{fig:2_holes_time_control})). The total number of time increments drops from 63 to 40.

As explained previously, the dissipated energy versus prescribed displacement curves (Figure (\ref{fig:dissi}))  obtained in the reference case 3 (very small prescribed arc-length) and in Test case 1 (coarse prescribed arc-length) are very different (incorrect equilibrium path in the second case). When using the time control strategy, $\nu_{d}^{time,dd}$ being successively set to $1 \times 10^{-2}$ and $5 \times 10^{-3}$, the correct equilibrium path is followed. In addition, the dissipated energy versus prescribed displacement curves gets closer to the one obtained in the reference simulation when the value of $\nu_{d}^{time,dd}$ decreases.

\subsection{Discussion}

The threshold value found numerically here can be compared to the one prescribed to ensure the convergence of the iterative resolution strategy at each computation time, $\nu^{iter}_d$. As explained in Section \ref{sec:latin}, the value $\nu^{iter}_d$ which ensures a sufficient convergence of the LaTIn solver can be obtained empirically by performing time independent benchmark tests (for instance the first time step of a delamination analysis). The time control strategy developed in this paper consists in monitoring the residual of the reference problem equations continuously during the time analysis, the measure used at any time being a time independent version of $\nu^{iter}$. Hence, it is not surprising to find out in the numerical examples that the higher value $\nu^{time,dd}_d$ permitting to follow the correct equilibrium path is the value of $\nu^{iter}$ which permits to obtain the convergence of the global informations (position of the crack fronts) at each computation time. Hence, applying the time control procedure only requires the prior knowledge of indicator $\nu^{iter}_d$.

